\documentclass[11pt,a4paper]{article}
\usepackage[cp1251]{inputenc}
\usepackage[english,russian]{babel}
\usepackage{oldlfont}

\usepackage{amsmath}
\usepackage{amsfonts}
\usepackage{amssymb}
\usepackage{array,longtable}
\usepackage{amscd}

\begin{document}

\begin{center}
\rule{0pt}{4cm}
    {\huge {\bf On spectrum and approximations \\
    of one class of sign-symmetric matrices}} \\[1cm]
\end{center}
\begin{center}
    {\Large Olga Y. Kushel \\
Department of mechanics and mathematics, \\
Belorussian State University, \\ Nezavisimosti sq., 4, 220050,
Minsk, Belarus, \\[0.15cm]
 e-mail: kushel@mail.ru}
\end{center}
\medskip

\begin{center}
{\Large May 2009}
\end{center}

\bigskip

\begin{center}
{\bf Abstract.}
\end{center}

A new class of sign-symmetric matrices is introduced in this
paper. Such matrices are named $\mathcal J$--sign-symmetric.
 The spectrum of a $\mathcal J$--sign-symmetric irreducible matrix
is studied under assumptions that its second compound matrix is
also $\mathcal J$--sign-symmetric and irreducible. The conditions,
when such matrices have complex eigenvalues on the largest
spectral circle, are given. The existence of two positive simple
eigenvalues $\lambda_1 > \lambda_2 > 0$ of a $\mathcal
J$--sign-symmetric irreducible matrix $A$ is proved under some
additional conditions. The question, when the approximation of a
$\mathcal J$--sign-symmetric matrix with a $\mathcal
J$--sign-symmetric second compound matrix by strictly $\mathcal
J$--sign-symmetric matrices with strictly $\mathcal
J$--sign-symmetric compound matrices is possible, is also studied
in this paper.

\medskip


\medskip


\medskip

\section{Introduction}

 First remind some well-known
definitions and statements of the theory of positive definite
matrices. The matrix $\mathbf{A}$ of a linear operator
$A:{\mathbb{R}}^n \rightarrow {\mathbb{R}}^n$ is called {\it
positive (strictly positive) definite}, if the inequality
$x^T{\mathbf A}x \geq 0$ (respectively $> 0$) is true for any
vector $x \in {\mathbb{R}}^n$. The following criterion of positive
(strictly positive) definiteness is known (see, for example, [1],
p. 47, theorem 8): {\it in order that a real symmetric matrix
$\mathbf{A}$ be positive (strictly positive) definite, it is
necessary and sufficient that all its principal minors be
nonnegative (respectively positive)}. The statement, that {\it all
the eigenvalues of a positive (strictly positive) definite matrix
are nonnegative (positive)}, is also well-known (see, for example,
[1], p. 46, theorem 7).

The following definitions, which generalize the concept of
positive (strictly positive) definiteness, were introduced in
paper [2] by M. Fiedler, V.P. Pt\'{a}k. The matrix $\mathbf{A}$ is
called {\it $P_0$--matrix} ($P$--matrix), if all its principal
minors of any order are nonnegative (respectively positive).

It's obvious, that any symmetric $P$--matrix is strictly positive
definite, and, as it follows, all its eigenvalues are positive.
However, the statement of reality, and all the more, of positivity
of the spectrum is not correct for an arbitrary (not necessarily
symmetric) $P$--matrix. (In this case one can only give the upper
bound of the argument of an arbitrary complex eigenvalue. Such a
statement was proved by R.B. Kellog in paper [3]). A natural
question arises, what additional conditions is it necessary to
impose on non-symmetric $P$--matrix in order to prove the
positivity of its spectrum.

The answer on this question was partly given in the 30th of the XX
century by F.R. Gantmacher and M.G. Krein. A class of strictly
totally positive matrices (i.e. matrices, {\it all} minors of
which are positive) was studied by this authors in monograph [1].
The following theorem was proved in [1] for the case of strictly
totally positive matrices.

{\bf Theorem A (Gantmacher--Krein)}. {\it If the matrix
$\mathbf{A}$ of a linear operator $A: {\mathbb{R}}^{n} \rightarrow
{\mathbb{R}}^{n}$ is strictly totally positive, then the operator
$A$ has $n$ positive simple eigenvalues $0 < \lambda_n < \ldots <
\lambda_2 < \lambda_1$, with the positive eigenvector $x_1$
corresponding to the maximal eigenvalue $\lambda_1$, and the
eigenvector $x_j$, which has exactly $j - 1$ changes of sign,
corresponding to the $j$-th eigenvalue $\lambda_j$.}

 The hypothesis, that it's not necessary to impose the positivity of {\it
all} minors for the positivity of all the eigenvalues of a
non-symmetric $P$--matrix, was made by D.M. Kotelyanskii. The
following theorem was proved in paper [4] by D.M. Kotelyanskii:
{\it if all the principal minors, and all the minors, obtained
from the principal minors by deleting one row and one column with
different numbers, of a matrix $\mathbf{A}$ (not necessarily
symmetric) are positive, then all the eigenvalues of $\mathbf{A}$
are simple and positive.}

So, the statement of the Gantmacher--Krein theorem of the
positivity and the simplicity of the eigenvalues was spread on the
wider class of matrices.

However, the statement, proved by Kotelyanskii, also doesn't cover
all the classes of matrices with positive spectrum. Remember the
following definitions for studying another classes of such
matrices. A $n \times n$ matrix $\mathbf{A}$ is called {\it
sign-symmetric}, if for any indices $i,j \in \{1, \ \ldots, \ n\}$
the inequality $a_{ij}a_{ji} \geq 0$ is true. (Such a definition
is used in papers by some modern authors (see, for example,
[5]--[6]), where such matrices are called "weakly
sign-symmetric"). It's obvious, that every nonnegative matrix is
sign-symmetric, i.e. the class of nonnegative matrices belongs to
the class of sign-symmetric matrices. Remind also the definition
of strict sign-symmetricity. A matrix $\mathbf{A}$ is called {\it
strictly sign-symmetric}, if for any indices $i,j \in \{1, \
\ldots, \ n\}$ the strict inequality $a_{ij}a_{ji} > 0$ is true.
In this case every positive matrix appears to be strictly
sign-symmetric.

Later on introduce the class of {\it totally sign-symmetric}
matrices. A matrix $\mathbf{A}$ is called {\it totally
sign-symmetric}, if for any $k = 1, \ \ldots, \ n$ and any sets
$\alpha = \{i_1, \ \ldots, \ i_k\}$ and $\beta = \{j_1, \ \ldots,
\ j_k \}$, for which $1 \leq i_1 < \ldots < i_k \leq n, \ 1 \leq
j_1 < \ldots < j_k \leq n$, the inequality $ A\begin{pmatrix}
  \alpha \\
  \beta
\end{pmatrix} A\begin{pmatrix}
  \beta  \\
  \alpha
\end{pmatrix} \geq 0$ is true, i.e. the product of
any two minors, which are symmetric with respect to the principal
diagonal, is non-negative. The class of totally positive matrices
belongs to the class of totally sign-symmetric matrices.

 In order to avoid a confusion, note, that sign-symmetricity of a matrix
was understood precisely as total sign-symmetricity in paper [7]
by D.M. Kotelyanskii, where the concept of sign-symmetricity was
first introduced, and also in papers by D. Hershkowitz, N. Keller,
O. Holtz and other authors (see., for example, [8]--[9]).

Introduce the class of {\it strictly totally sign-symmetric}
matrices.
 A matrix $\mathbf{A}$ is called {\it strictly totally sign-symmetric},
if $\mathbf{A}$ does not contain any zero elements, and for any
sets $\alpha = \{i_1, \ \ldots, \ i_k\}$ and $\beta = \{j_1, \
\ldots, \ j_k \}$, for which $1 \leq i_1 < \ldots < i_k \leq n, \
1 \leq j_1 < \ldots < j_k \leq n$, the inequality $
A\begin{pmatrix}
  \alpha \\
  \beta
\end{pmatrix} A\begin{pmatrix}
  \beta  \\
  \alpha
\end{pmatrix} > 0$ is true, i.e. the matrix $\mathbf{A}$ does not
have zero minors, and the product of any two minors, which are
symmetric with respect of the principal diagonal, is positive. The
class of strictly totally positive matrices belongs to the class
of strictly totally sign-symmetric matrices.

 It's easy to see, that neither sign-symmetricity,
no total sign-symmetricity of a $P$--matrix does not guarantee the
reality of its spectrum. However, some authors describes spectral
properties of totally sign-symmetric $P$--matrices, in particular,
D. Carlson in paper [10] proved the statement, that {\it totally
sign-symmetric $P$-matrices are positive stable, i.e., all its
eigenvalues have positive real parts}.

 A problem of the separation of a new subclass of
matrices with positive spectra in the class of strictly totally
sign-symmetric $P$--matrices is studied in this paper. It's also
supposed to separate a subclass in the class of totally
sign-symmetric $P_0$--matrix, and to study, when a matrix, which
belongs to the separated subclass, has complex eigenvalues, and
when all its eigenvalues are real. It's also supposed to study,
when a totally sign-symmetric $P_0$-matrix, which belong to the
separated subclass, will not be positive stable, and to describe
its spectrum.

Besides, the question, when the approximation of  totally
sign-symmetric matrices by strictly totally sign-symmetric
matrices, is possible, is still open (see [9]). It is supposed to
give the answer on this question for totally sign-symmetric
$P_0$--matrices, which belongs to the separated subclass. Later on
we'll be restricted to the case, when the conditions will be
imposed only on the elements of the matrix and on its minors of
the second order.

\section{${\mathbb R}^n$ as the space of functions defined on a finite
support. Tensor and exterior powers of ${\mathbb R}^n$.} Later on,
as the author believes, it will be more convenient to consider the
space ${\mathbb R}^n$ as the space of functions, defined on the
finite set $\{1, \ \ldots, \ n \}$.

 Let us examine the set of the indices $\{1, \ \ldots, \ n\}$ and
the space ${\mathbb X}$ of functions $x: \{1, \ \ldots, \ n\}
\rightarrow {\mathbb R}$, defined on it. It is obvious, that the
space ${\mathbb X}$ is isomorphic to ${\mathbb R}^n$. The basis in
the space ${\mathbb X}$ consists of the functions $e_i$, for which
$e_i(j) = \delta_{ij}$.

The tensor square ${\mathbb X} \otimes {\mathbb X}$ of the space
${\mathbb X}$ is the space of all functions, defined on the set
$\{1, \ \ldots, \ n\} \times \{1, \ \ldots, \ n\}$, which consists
of $n^2$ pairs of the form $(i,j)$, where $i, j \in \{1, \ \ldots,
\ n\}$. The tensor product of functions $x$ and $y$ from ${\mathbb
X}$ is defined as the function $x \otimes y$, which acts according
to the rule:
$$(x \otimes y)(i,j) = x(i)y(j). $$

It's known (see, for example, [11], [12]), that all the possible
tensor products $e_i \otimes e_j$, $i,j = 1, \ldots, n$ of the
initial basic functions forms a basis in the space ${\mathbb
X}\otimes{\mathbb X}$. In this case the basic functions $e_i
\otimes e_j, \ \ i,j = 1, \ \ldots, \ n$ are numerated in order,
corresponding to the lexicographic numeration of the pairs
$(i,j)$. As it follows, the tensor square of the set ${\mathbb
R}^n$ can be considered as the space ${\mathbb R}^{n^2}$.

 Examine the exterior square ${\mathbb X}\wedge {\mathbb X}$
of the space ${\mathbb X}$. The exterior square ${\mathbb X}\wedge
{\mathbb X}$ is a subspace of the space ${\mathbb X} \otimes
{\mathbb X}$, which consists of all antisymmetric functions (i.e.
functions $f(i,j)$, for which the equality $f(i,j) = - f(j,i)$ is
true), defined on the set $\{1, \ \ldots, \ n\} \times \{1, \
\ldots, \ n\}$. It is known, that the exterior square ${\mathbb
X}\wedge {\mathbb X}$ coincides with the linear span of all
exterior products $x \wedge y \ \ (x,y \in {\mathbb X})$, which
acts according to the rule:
$$(x \wedge y)(i,j) = (x \otimes y)(i,j) - (y \otimes x)(i,j) =
x(i)y(j) - x(j)y(i).$$

Let $W$ be a subset of $\{1, \ \ldots, \ n\} \times \{1, \ \ldots,
\ n\}$, which satisfies the following conditions:$$ W\cup
\widetilde{W} = (\{1, \ \ldots, \ n\} \times \{1, \ \ldots, \
n\}); \eqno(1)$$
$$W \cap \widetilde{W} = \Delta. \eqno(2)$$ (Here $\widetilde{W}
= \{(j,i): \ (i,j) \in W\}$; $\Delta = \{(i,i): \ i =1, \ \ldots,
\ n\}$).

The following statement about the exterior square of the space
${\mathbb X}$ is true: {\it the space ${\mathbb X}\wedge {\mathbb
X}$ of antisymmetric functions, defined on the set $\{1, \ \ldots,
\ n\} \times \{1, \ \ldots, \ n\}$, is isomorphic to the space
${\mathbb X} (W\setminus\Delta)$ of real-valued functions, defined
on the set $W\setminus\Delta$.}

 Really, any function, defined on the set
$W\setminus\Delta$, where $W$ is a subset of $\{1, \ \ldots, \ n\}
\times \{1, \ \ldots, \ n\}$, which satisfies conditions (1) and
(2), can be extended as an antisymmetric function to $ \{1, \
\ldots, \ n\} \times \{1, \ \ldots, \ n\}$ by the unique way. Then
the received antisymmetric function is supposed to be equal to
zero on the set $\Delta = \{(i,i): \ i =1, \ \ldots, \ n\}$.

It's obvious, that the set $W$ is not uniquely defined. However,
its power remains constant and can be calculated by the following
way. The equality $$ N(W\cup \widetilde{W}) = N(\{1, \ \ldots, \
n\} \times \{1, \ \ldots, \ n\}) = n^2$$ follows from condition
(1). The equality $$N(W\cap \widetilde{W}) = N(\Delta) = n$$
follows from condition (2). The equality $N(W) = N(\widetilde{W})$
follows from the definition of the set $\widetilde{W}$. Then,
taking into account, that $N(W\cup \widetilde{W}) = N(W) +
N(\widetilde{W}) - N(W\cap \widetilde{W}) = 2N(W) - N(\Delta)$,
we'll get the equality:
$$N(W) = \frac{N(W\cup \widetilde{W}) + N(\Delta)}{2} = \frac{n^2
+ n}{2}.$$

In turn, the following equality is true for the power of the set
$W\setminus\Delta$:

 $$N(W\setminus\Delta) = N(W) - N(\Delta) = \frac{n^2 - n}{2} = C_n^2.$$

The following sequence of isomorphisms comes out from the above
reasoning:

$${\mathbb R}^n \wedge {\mathbb R}^n = {\mathbb X}(W \setminus
\Delta) = {\mathbb R}^{C_n^2}.$$

\section{ The set $W$ and binary relations on the set of indices $\{1,
\ldots, n \}$.} It's known, that the separation of a subset in a
Cartesian square of a set can be considered as a definition of a
binary relation on the initial set, the inverse statement is also
true (see, for example, [13]). As it follows, the definition of a
subset $W$ in the set $\{1, \ \ldots, \ n\} \times \{1, \ \ldots,
\ n\}$ is equivalent to the definition of a binary relation on the
set of indices $\{1, \ \ldots, \ n\}$. Describe the properties of
the set $W$ in terms of the corresponding binary relation. Later
on we shall use the following definitions.

 Consider the indices $i,j \in \{1, \ \ldots, \ n\}$ satisfy the
binary relation $W$ and note $ i \stackrel{W}{\prec} j$ if and
only if the pair $(i,j)$ belongs to the set $W$. In this case the
inverse relation is defined by the set $\widetilde{W}$.

 A binary relation $\stackrel{W}{\prec}$ on the set of indices $\{1, \ \ldots,
\ n\}$ is called:

- {\it reflexive}, if $i \stackrel{W}{\prec} i$ for any $i \in
\{1, \ \ldots, \ n\}$. It means, that the set $\Delta = \{(i,i): \
i = 1, \ \ldots, \ n\}$ belongs to both $W$ and $\widetilde{W}$.

- {\it antisymmetric}, if the equality $i = j$ follows from $i
\stackrel{W}{\prec} j, \ j \stackrel{W}{\prec} i$ for any $i, j
\in \{1, \ \ldots, \ n\}$. It means, that $W\cap \widetilde{W} =
\Delta$.

- {\it transitive}, if $i \stackrel{W}{\prec} k$ follows from $i
\stackrel{W}{\prec} j$ and $j \stackrel{W}{\prec} k$ for any $i,
j, k \in \{1, \ \ldots, \ n\}$. It means, that the inclusion
$(i,k) \in W$ follows from the inclusions $(i,j) \in W$ and $(j,k)
\in W$ for any $i, j, k \in \{1, \ \ldots, \ n\}$.

- {\it connected}, if either $i \stackrel{W}{\prec} j$ or $j
\stackrel{W}{\prec} i$ is true for any pair $(i,j)$. It means,
that $W\cup \widetilde{W} = \{1, \ \ldots, \ n\} \times \{1, \
\ldots, \ n\}$.

If a binary relation $\stackrel{W}{\prec}$ is reflexive,
antisymmetric, transitive and connected, then it is called a
linear order relation (see, for example, [14]).

The following statement comes out from the above reasoning.

{\it Let the set $W \subset \{1, \ \ldots, \ n\} \times \{1, \
\ldots, \ n\}$ satisfy conditions (1) and (2). Then it defines a
connected antisymmetric reflexive binary relation on the set $\{1,
\ \ldots, \ n\}$. If, in addition, the inclusion $(i,k) \in W$
follows from the inclusions $(i,j) \in W$ and $(j,k) \in W$ for
any $i, j, k \in \{1, \ \ldots, \ n\}$, then the relation, defined
by the set $W$, is a linear order relation.}

Show that the inverse statement is also true. Let a connected
antisymmetric reflexive binary relation $\stackrel{W}{\prec}$ be
defined on the set $\{1, \ \ldots, \ n\}$. Then we can define $W$
and $ \widetilde{W}$ by the following way:
$$W = \{(i, j) \in \{1, \ \ldots, \ n\} \times \{1, \ \ldots, \
n\} : i \stackrel{W}{\prec} j\}; $$ $$\widetilde{W} = \{(i, j) \in
\{1, \ \ldots, \ n\} \times \{1, \ \ldots, \ n\}: j
\stackrel{W}{\prec} i\}. $$ Properties (1) and (2) of the sets $W$
and $\widetilde{W}$ follow from the properties of the relation
$\stackrel{W}{\prec}$.

 The set $M = \{(i,j) \in \{1, \ \ldots, \ n\} \times \{1, \ \ldots, \
n\}: \ i \leq j \}$ is considered as $W$ in the classical theory
of totally positive and oscillatory matrices (see, for example,
[1]). It corresponds the natural linear order relation on $\{1, \
\ldots, \ n\}$.

It's easy to see, that the number of all the possible methods of
constructing the set $W$, which satisfies conditions (1) and (2),
is equal to $2^{C_n^2}$. Respectively, we'll get $2^{C_n^2}$ ways
of the realization of the exterior square ${\mathbb X} \wedge
{\mathbb X} = {\mathbb X}(W \setminus \Delta)$.

Note, that the property of the set $W$, which defines the
transitivity of the corresponding binary relation, is not
necessary for equalities (1) and (2) to be true. However, the sets
$W$, which define linear order relations on $\{1, \ \ldots, \
n\}$, will play an important role in what follows. The number of
such sets is equal to the number of all the possible permutations
of the set $\{1, \ \ldots, \ n\}$, i.e. $n!$.

\section{Basis in the exterior power of ${\mathbb R}^{n}$ and its connection
with the construction of the set $W$.} The following statement is
well-known (see, for example, [11], [12]): {\it if $e_1, \ \ldots,
\ e_n$ is a basis in ${\mathbb{R}}^{n}$, then all the possible
exterior products of the form $e_i \wedge e_j$, where $1 \leq i <
j \leq n $, forms a basis in the exterior square ${\mathbb{R}}^{n}
\wedge {\mathbb{R}}^{n}$ of the space $ {\mathbb{R}}^{n}$.}

However, note, that there exist other bases, which consist of
exterior products of the initial basic vectors, in the space
${\mathbb{R}}^{n} \wedge {\mathbb{R}}^{n}$, besides the canonical
basis $\{e_i \wedge e_j\}$, where $1 \leq i < j \leq n $. Such
bases are constructed by the following way: one arbitrary element
is selected from every pair of the opposite elements $e_i \wedge
e_j$ and $e_j \wedge e_i \ \ (i \neq j)$. So, one can construct
$2^{C_n^2}$ different bases.

Let us prove the following lemma about the connection between a
set $W$ and a basis in ${\mathbb R}^n \wedge {\mathbb R}^n$.

{\bf Lemma 1.} {\it Every set $W$, which satisfies conditions (1)
and (2), uniquely defines a basis in ${\mathbb R}^n \wedge
{\mathbb R}^n$, which consists of the exterior products of the
initial basic vectors. The inverse statement is also true: every
basis in ${\mathbb R}^n \wedge {\mathbb R}^n$, which consists of
the exterior products of the initial basic vectors, uniquely
defines a subset $W$ in $\{1, \ \ldots, \ n\}\times \{1, \ \ldots,
\ n\}$, which satisfies conditions (1) and (2)}.

{\bf Proof.} First prove, that every set $W$, which satisfies
conditions (1) and (2), uniquely defines a basis in ${\mathbb X}
\wedge {\mathbb X}$. Examine the system $\Lambda$, which consists
of the exterior products $e_i \wedge e_j, \ \ (i,j)\in W\setminus
\Delta$. Show, that $\Lambda$ is a basis in ${\mathbb X} \wedge
{\mathbb X}$. It's enough for this to show, that the restrictions
of the functions from $\Lambda$ to the set $W \setminus \Delta$
form a basis in ${\mathbb X}(W \setminus \Delta)$. Really, examine
the value of an arbitrary function $e_i \wedge e_j \in \Lambda$ on
an arbitrary pair of indices $ (k,l) \in W \setminus \Delta$. It's
easy to see, that
$$(e_i \wedge e_j)(k,l)
 = \left\{\begin{array}{cc} 1, &
\mbox{if $(i,j) = (k,l)$ };
\\[10pt] 0, &  \mbox{otherwise.}\end{array}\right.$$
As it follows, the functions from $\Lambda$ are linearly
independent, and since the system $\Lambda$ contains $C_n^2$
functions, this system is complete.

 Prove the inverse statement. Given a basis $\Lambda$ of the space ${\mathbb X} \wedge {\mathbb
X}$, which consists of the exterior products of the initial basic
functions. Construct the set $W$ by the following way: $(i,j) \in
W$, if and only if either $i = j$ or $e_i \wedge e_j \in \Lambda$.
Show, that such a set satisfies conditions (1) and (2). First
verify condition (1). Assume, that there exists a pair $(i_0, j_0)
\ \ i_0 \neq j_0$, which belongs to $W \cap \widetilde{W}$. In
this case it follows from the definition of the set
$\widetilde{W}$, that the pair $(j_0, i_0)$ also belongs to $W
\cap \widetilde{W}$, i.e. both the pairs $(i_0, j_0)$ and $(j_0,
i_0)$ belongs to the set $W$. As it follows, both $e_{i_0} \wedge
e_{j_0}$, and $e_{j_0} \wedge e_{i_0}$ belongs to the linearly
independent system $\Lambda$. We came to the contradiction, i.e.
$e_{i_0} \wedge e_{j_0} = - (e_{j_0} \wedge e_{i_0})$, and, as it
follows, this functions are linearly dependent. Verify condition
(2). Let a pair $(i_0,j_0) \ \ i_0 \neq j_0$, which belongs to $
\{1, \ \ldots, \ n\} \times \{1, \ \ldots, \ n\}$, but does not
belong to $W \cup \widetilde{W}$, does exist. Then the pair $(j_0,
i_0)$ also does not belong to $W \cup \widetilde{W}$. As it
follows, neither the function $e_{i_0} \wedge e_{j_0}$, no the
function $e_{j_0} \wedge e_{i_0}$ does not belong to the system
$\Lambda$. Let us add the function $e_{i_0} \wedge e_{j_0}$ to the
system $\Lambda$. It's easy to see, that the obtained system is
linearly independent. This fact contradicts the condition, that
$\Lambda$ is the maximal linearly independent system in the space
${\mathbb X} \wedge {\mathbb X}$. $\square$

Later on we'll call the basis $\{e_i \wedge e_j\}_{(i,j)\in W
\setminus \Delta}$, constructed with respect to the set $W$, {\it
a $W$--basis}. Note, that the elements of a $W$--basis are
numerated in the lexicographic order. Let us give an example of a
$W$--basis.

{\bf Example 1.} Let $M = \{(i,j) \in \{1, \ \ldots, \ n\} \times
\{1, \ \ldots, \ n\}: \ i \leq j\}$. Then \linebreak $M\setminus
\Delta = \{(i,j) \in \{1, \ \ldots, \ n\} \times \{1, \ \ldots, \
n\}: \ i < j\}$, and the corresponding $W$--basis is a set of
those exterior products $e_i \wedge e_j$ of the initial basic
vectors, for which $i < j$. Such a $W$--basis is a canonical basis
in the space ${\mathbb R}^n \wedge {\mathbb R}^n$. It is studied
in many papers (see, for example, [1], [12]).

\section{Exterior power of a linear operator in ${\mathbb R}^n$
and its matrix in the $W$--basis}
 The exterior square $A \wedge A$ of the operator
$A: {\mathbb X}\rightarrow {\mathbb X}$ acts in the space
${\mathbb X}\wedge{\mathbb X}$ according to the rule: $$ (A \wedge
A)(x \wedge y) = Ax \wedge Ay. $$ Later on we shall study spectral
properties of the operator $A$, and we shall require its exterior
square $A \wedge A$ to leave invariant a cone in ${\mathbb
X}\wedge{\mathbb X}$ (in this case the spectral radius $\rho(A
\wedge A)$ of the operator $A \wedge A$ is an eigenvalue of $A
\wedge A$). It's enough for this the matrix of the operator $A
\wedge A$ to be positive in some $W$--basis in ${\mathbb
X}\wedge{\mathbb X}$. Then the operator $A \wedge A$ leaves
invariant the cone, spanned on the vectors of this $W$--basis.

 Let the operator $A$ be defined by the matrix ${\mathbf A} = \{a_{ij}\}_{i,j =
1}^n$ in the basis $\{e_i\}_{i = 1}^n$. Let us study the matrix of
the operator $A \wedge A$ in a $W$--basis, constructed with
respect to a set $W$, which satisfies conditions (1) and (2).
First remember the following definition.

A minor $A\begin{pmatrix}
  i & j \\
  k & l
\end{pmatrix}$, formed of the rows with numbers $i$ and $j$
and the columns with numbers $k$ and $l$, where indices $i, \ j, \
k, \ l$ take any values from $\{1, \ldots, n\}$, is called {\it a
generalized minor of the second order of the matrix ${\mathbf
A}$}. Note, that if we change places of the rows (columns), the
sign of the minor will be changed, that is why when $i = j \ $ ($k
= l$), the minor $A\begin{pmatrix}
  i & j \\
  k & l
\end{pmatrix}$ is equal to zero.

Call the matrix, which consists of generalized minors of the
second order $A\begin{pmatrix}
  i & j \\
  k & l
\end{pmatrix}$, where $ (i,j), \ (k,l) \in (W \setminus \Delta)$, numerated in the lexicographic order, {\it a
$W$--matrix}, and denote it ${\mathbf A}_W^{(2)}$.

{\bf Example 2.} Let $W = M = \{(i,j) \in \{1, \ \ldots, \ n\}
\times \{1, \ \ldots, \ n\}: \ i \leq j\}$. Then the corresponding
$W$--matrix is a matrix, which consists of minors
$A\begin{pmatrix}
  i & j \\
  k & l
\end{pmatrix}$, where $i < j, \ k < l$, i.e.
{\it the second compound matrix}.

Let us prove the following theorem illustrating the connection
between a $W$--matrix and a matrix of the exterior square of the
operator $A$.

{\bf Theorem 1}. \textit{Let the operator $A$ be defined by the
matrix ${\mathbf A} = \{a_{ij}\}_{i,j = 1}^n$ in the basis $e_{1},
\ldots, e_{n}$. Let a subset $W \subset \{1, \ \ldots, \ n\}
\times \{1, \ \ldots, \ n\}$ satisfy conditions (1) and (2). Then
the matrix of the exterior square $A \wedge A$ of the operator $A$
in the $W$--basis $\{e_i \wedge e_j\}_{(i,j)\in W \setminus
\Delta}$ coincides with the $W$--matrix $ {\mathbf A}_W^{(2)}$.}

{\bf Proof.} Let us compare the columns of the matrices, using the
fact, that $A(e_{k}) = \sum\limits_{i = 1}^{n}a_{ik}e_{i}$ for $k
= 1, \ \ldots, \ n$. Examine the column of the matrix of the
operator $A \wedge A$ with an arbitrary number $\alpha$. The
number $\alpha$ in the lexicographic numeration corresponds to a
pair of indices $(i, j) \in W \setminus \Delta$. Prove, that this
column coincides with the column of the $W$--matrix ${\mathbf
A}_W^{(2)}$, which has the same number:
$$(A\wedge A)(e_{i}\wedge e_{j}) = Ae_{i}\wedge Ae_{j} =
\left(\sum_{k =
 1}^{n}a_{ki}e_{k}\right)\wedge \left(\sum_{l = 1}^{n}a_{lj}e_{l}\right) =
\sum_{k,l = 1}^{n}a_{ki}a_{lj}(e_{k}\wedge  e_{l}) = $$ $$ =
\sum_{(k,l) \in (W\setminus\Delta)}a_{ki}a_{lj}(e_{k}\wedge e_{l})
+ \sum_{k = l = 1}^{n}a_{ki}a_{lj}(e_{k}\wedge  e_{l}) + \sum_{
(k,l) \in (\widetilde{W}\setminus\Delta)}a_{ki}a_{lj}(e_{k}\wedge
e_{l}) = $$ $$= \sum_{ (k,l) \in
(W\setminus\Delta)}a_{ki}a_{lj}(e_{k}\wedge e_{l}) + 0 - \sum_{
(k,l) \in (\widetilde{W}\setminus\Delta)}a_{ki}a_{lj}(e_{l}\wedge
e_{k}).$$ Change places of the indices $l$ and $k$ in the third
sum:
$$\sum_{(k,l) \in (W\setminus\Delta)}a_{ki}a_{lj}(e_{k}\wedge
e_{l}) - \sum_{(k,l) \in
(W\setminus\Delta)}a_{li}a_{kj}(e_{k}\wedge e_{l}) = \sum_{ (k,l)
\in (W\setminus\Delta)}(a_{ki}a_{lj} - a_{li}a_{kj})(e_{k}\wedge
e_{l}) =$$ $$= \sum_{(k,l) \in (W\setminus\Delta)}A\begin{pmatrix}
k & l \\ i & j \
\end{pmatrix} \ (e_{k}\wedge  e_{l}),$$ where $A\begin{pmatrix} k &
l \\ i & j \  \end{pmatrix}$ are the elements of the column with
the number $\alpha$ of the matrix ${\mathbf A}_W^{(2)}$. That is,
the matrix of the exterior square of the operator $A$ coincides
with the $W$--matrix ${\mathbf A}_W^{(2)}$.

{\bf Corollary.} {\it The matrix of the exterior square of the
operator $A$ in the basis $\{e_i \wedge e_j\}_{i < j}$ coincides
with the second compound matrix of the matrix $A$.}

 Let us prove the following theorem about the eigenvalues of a $W$--matrix.

{\bf Theorem 2.} \textit{Let $W$ be a set, which satisfies
conditions (1) and (2). Let $\{\lambda_{i}\}_{i = 1}^n$ be the set
of all eigenvalues of the matrix ${\mathbf A}$, repeated according
to multiplicity. Then all the possible products of the type
$\{\lambda_{i} \lambda_{j} \}$, where $1 \leq i < j \leq n$, forms
the set of all the possible eigenvalues of the corresponding
$W$--matrix $ {\mathbf A}_W^{(2)}$, repeated according to
multiplicity}.

{\bf Proof.} It follows from theorem 1, that the $W$--matrix
${\mathbf A}_W^{(2)}$ coincides with the matrix of the exterior
square of the operator $A$ in the corresponding $W$--basis for any
$W$, which satisfies conditions (1) and (2). The following
statement is true for the exterior square $A \wedge A$ of the
operator $A$ (see, for example, [12]): all the possible products
of the type $\{\lambda_{i} \lambda_{j} \}$, where $1 \leq i < j
\leq n$, forms the set of all the possible eigenvalues of $A
\wedge A $, repeated according to multiplicity. $\square$

Note, that in the case $W = M$ theorem 2 turns into the Kronecker
theorem (see [1], p. 80, theorem 23) about the eigenvalues of the
second compound matrix. The proof of the Kronecker theorem without
using exterior products is given in monograph [1].

\section{Classes of matrices. Basic definitions and statements}

The proof of the theorem А (Gantmacher--Krein) is based on the
wide-known result of Perron and Frobenius about the existence of
the largest positive eigenvalue of a non-negative irreducible
matrix $\mathbf{A}$. However, it's easy to see, that this result
is also correct for any matrix, similar to the matrix $\mathbf{A}$
(because of the spectrum of a matrix does not change when passing
to another basis). Hear a natural question arises: how can we see
if an arbitrary matrix is similar to some nonnegative matrix? Let
us formulate and prove a sufficient criterion of similarity, which
will be used later in this paper.

 First remember some definitions and statements from the matrix
theory (see, for example, [15]). A matrix ${\mathbf A}$ is called
{\it non-negative} ({\it positive}), if all its elements $a_{ij}$
are nonnegative (positive). The following statement (Perron's
theorem) is true for positive matrices: {\it let the matrix
$\mathbf A$ of a linear operator $A:{\mathbb{R}}^{n} \rightarrow
{\mathbb{R}}^{n}$ be positive. Then the spectral radius $\rho(A) >
0$ is a simple positive eigenvalue of the operator $A$, different
in modulus from the other eigenvalues. Besides, the eigenvector
$x_1$, corresponding to the eigenvalue $\lambda_1 = \rho(A)$, is
positive.}

Let us give the following definition, which generalizes the
definition of matrix positivity. A matrix $\mathbf A$ is called
{\it strictly ${\mathcal J}$--sign-symmetric}, if ${\mathbf A}$
does not contain zero elements and there exists such a subset
${\mathcal J} \subseteq \{1, \ \ldots, \ n\}$, that the inequality
$a_{ij} < 0$ is true if and only if one of the numbers $i$, $j$
belongs to the set ${\mathcal J}$, and the other belongs to the
set $\{1, \ \ldots, \ n\}\setminus {\mathcal J}$. It's obvious,
that positive matrices belong to the class of strictly ${\mathcal
J}$--sign-symmetric matrices. Note, that the set ${\mathcal J}$ in
the definition of strict ${\mathcal J}$--sign-symmetricity is
uniquely defined (up to the set $\{1, \ \ldots, \ n\}\setminus
{\mathcal J}$).

Let us prove the following statement for ${\mathcal
J}$--sign-symmetric matrices.

{\bf Theorem 3.} {\it Let ${\mathbf A}$ be a strictly ${\mathcal
J}$--sign-symmetric matrix.
 Then it can be represented in the following form: $${\mathbf A} = {\mathbf
D} \widetilde{{\mathbf A}} {\mathbf D}^{-1},$$ where
$\widetilde{{\mathbf A}}$ is a positive matrix, ${\mathbf D}$ is a
diagonal matrix, which diagonal elements are equal to $\pm 1$.}

{\bf Proof.} Let ${\mathbf A}$ be a strictly ${\mathcal
J}$--sign-symmetric matrix. Then there exists such a subset
${\mathcal J} \subseteq \{1, \ \ldots, \ n\}$, that the inequality
$a_{ij} < 0$ is true if and only if one of the numbers $i$, $j$
belongs to the set ${\mathcal J}$, and the other belongs to the
set $\{1, \ \ldots, \ n\}\setminus {\mathcal J}$. In the case of
${\mathcal J} = \{1, \ \ldots, \ n\}$ or ${\mathcal J} =
\emptyset$, it's easy to see, that all the elements of the matrix
${\mathbf A}$ are positive. Then $\widetilde{{\mathbf A}}$
coincides with the matrix ${\mathbf A}$, and ${\mathbf D}$ is an
identity matrix.

Let ${\mathcal J} \neq \{1, \ \ldots, \ n\}$ and ${\mathcal J}
\neq \emptyset$. Define the diagonal matrix ${\mathbf D}$ by the
following way:

$$d_{ii}
 = \left\{\begin{array}{cc} - 1, &
\mbox{if $i \in {\mathcal J}$ };
\\[10pt] 1, &  \mbox{otherwise.}\end{array}\right.$$

In this case it's obvious, that ${\mathbf D}^{-1} = {\mathbf D}$.

It's easy to see, that if $i \in \{1, \ \ldots, \ n\}\setminus
{\mathcal J}$, then the $i$-th column of the product ${\mathbf
A}{\mathbf D}$ coincides with the $i$-th column of the initial
matrix ${\mathbf A}$, and if $i \in {\mathcal J}$, then the $i$-th
column of ${\mathbf A}{\mathbf D}$ coincides with the $i$-th
column of ${\mathbf A}$, taken with the opposite sign. In turn, if
$i \in \{1, \ \ldots, \ n\}\setminus {\mathcal J}$, then the
$i$-th row of ${\mathbf D}{\mathbf A}$ is equal to the $i$-th row
of ${\mathbf A}$, and if $i \in {\mathcal J}$ then the $i$-th row
of ${\mathbf D}{\mathbf A}$ is equal to the $i$-th row of
${\mathbf A}$, taken with the opposite sign.

Show, that the matrix $\widetilde{{\mathbf A}} = {\mathbf
D}^{-1}{\mathbf A} {\mathbf D}$, is positive. Examine its
arbitrary element $\widetilde{a}_{ij}$. One of the following three
cases takes place:
\begin{enumerate}

\item[\rm 1.] Both the indices $i,j$ belong to the set
$\{1, \ \ldots, \ n\}\setminus{\mathcal J}$. In this case it's
easy to see, that the element $\widetilde{a}_{ij}$ of the matrix
$\widetilde{{\mathbf A}}$ is equal to the corresponding element
$a_{ij}$ of the strictly ${\mathcal J}$--sign-symmetric matrix
${\mathbf A}$. It follows from the definition of strict ${\mathcal
J}$--sign-symmetricity, that the element $a_{ij}$ is positive.

\item[\rm 2.] One of the indices $i,j$
belongs to the set ${\mathcal J}$, and the other belongs to the
set $\{1, \ \ldots, \ n\}\setminus{\mathcal J}$. Let us take $i
\in {\mathcal J}$, $j \in \{1, \ \ldots, \ n\}\setminus{\mathcal
J}$ (the second case is studied by analogy). Then the element
$a_{ij}$ of the matrix ${\mathbf A}$ changes its sign (with the
$i$-th row), when multiplying from the left by the matrix
${\mathbf D}^{-1} = {\mathbf D}$, and it remains the same (with
the $j$-th column) when multiplying from the right by ${\mathbf
D}$. As it follows, the element $\widetilde{a}_{ij}$ of the matrix
$\widetilde{{\mathbf A}}$ is equal to the element $a_{ij}$ of the
strictly ${\mathcal J}$--sign-symmetric matrix ${\mathbf A}$,
taken with the opposite sign. It follows from the definition of
strict ${\mathcal J}$--sign-symmetricity, that if one of the
indices $i,j$ belongs to the set ${\mathcal J}$, and the other
belongs to the set $\{1, \ \ldots, \ n\}\setminus{\mathcal J}$,
then the element $a_{ij}$ is negative. As it follows,
$\widetilde{a}_{ij} = - a_{ij}$ is positive.

\item[\rm 3.] Both the indices $i,j$ belong to the set
${\mathcal J}$. In this case the element $a_{ij}$ changes its sign
twice: with the $i$-th row when multiplying from the left by
$\mathbf D$, and with the $j$-th column when multiplying from the
right by $\mathbf D$. As it follows, $\widetilde{a}_{ij} =
a_{ij}$, and, using the definition of strict ${\mathcal
J}$--sign-symmetricity, we'll get, that $a_{ij}$ is positive.
\end{enumerate} $\square$

{\bf Corollary.} {\it Let the matrix $\mathbf A$ of a linear
operator $A:{\mathbb{R}}^{n} \rightarrow {\mathbb{R}}^{n}$ be
strictly ${\mathcal J}$--sign-symmetric. Then the spectral radius
$\rho(A)> 0$ is a simple positive eigenvalue of the operator $A$,
different in modulus from the other eigenvalues.}

The class of positive matrices belongs to the more general class
of irreducible nonnegative matrices. To describe this class,
remember the following definition. A matrix ${\mathbf A}$ is
called {\it reducible}, if there exists a permutation of
coordinates such that: $${\mathbf P}^{-1}{\mathbf A}{\mathbf P} =
\begin{pmatrix}
  {\mathbf A}_1 & 0 \\
  {\mathbf B} & {\mathbf A}_2
\end{pmatrix}, \eqno(3)$$ where ${\mathbf P}$ is an $n \times n$ permutation matrix
(each row and each column have exactly one 1 entry and all others
0), ${\mathbf A}_1$, ${\mathbf A}_2$ are square matrices.
Otherwise the matrix $\mathbf A$ is called irreducible. The
Frobenius theorem, which generalizes the Perron theorem to the
case of nonnegative irreducible matrices, is widely known: {\it
let the matrix $\mathbf A$ of a linear operator $A$ be nonnegative
and irreducible. Then the spectral radius $\rho(A)
> 0$ is a simple positive eigenvalue of the operator $A$,
with the corresponding positive eigenvector $x_1$. If $h$ is a
number of the eigenvalues of the operator $A$, which are equal in
modulus to $\rho(A)$, then all of them are simple and they
coincide with the $h$-th roots of $(\rho(A))^h$. More than that,
the spectrum of the operator $A$ is invariant under rotations by
$\frac{2\pi}{h}$ about the origin.}

 The number $h$ of the eigenvalues, which are equal in modulus to
$\rho(A)$, is called {\it the index of imprimitivity} of the
irreducible operator $A$. The operator $A$ is called {\it
primitive}, if $h(A) = 1$, and {\it imprimitive}, if $h(A)
> 1$.

Let us generalize the definition of strict ${\mathcal
J}$--sign-symmetricity in order to receive a class of matrices,
which includes all non-negative matrices.

 A matrix ${\mathbf A}$ of a linear operator $A:
{\mathbb R}^n \rightarrow {\mathbb R}^n$ is called {\it ${\mathcal
J}$--sign-symmetric}, if there exists such a subset ${\mathcal J}
\subseteq \{1, \ \ldots, \ n\}$, that the inequality $a_{ij} \leq
0$ is true for any two numbers $i,j$, one of which belongs to the
set ${\mathcal J}$, and the other belongs to the set $\{1, \
\ldots, \ n\}\setminus {\mathcal J}$; and the strict inequality
$a_{ij} < 0$ is true only if one of the numbers $i$, $j$ belongs
to  ${\mathcal J}$, and the other belongs to $\{1, \ \ldots, \
n\}\setminus {\mathcal J}$. It's not difficult to see, that if the
matrix ${\mathbf A}$ is ${\mathcal J}$--sign-symmetric and
irreducible, then the set ${\mathcal J}$ in the definition of
${\mathcal J}$--sign-symmetricity is uniquely defined (up to the
set $\{1, \ \ldots, \ n\}\setminus {\mathcal J}$). If ${\mathbf
A}$ is reducible, then there is a finite number $k > 2$ of
possible ways of constructing the set ${\mathcal J}$.

Let us generalize theorem 3 to the case of ${\mathcal
J}$--sign-symmetric irreducible matrices.

{\bf Theorem 4.} {\it Let ${\mathbf A}$ be a ${\mathcal
J}$--sign-symmetric matrix. Then it can be represented in the
following form:
$${\mathbf A} = {\mathbf D} \widetilde{{\mathbf A}} {\mathbf
D}^{-1},$$ where $\widetilde{{\mathbf A}}$ is a nonnegative
matrix, ${\mathbf D}$ is a diagonal matrix, which diagonal
elements are equal to $\pm 1$. More than that, if ${\mathbf A}$ is
irreducible, then $\widetilde{{\mathbf A}}$ is also irreducible.}

{\bf Proof}. The proof of the fact, that a ${\mathcal
J}$--sign-symmetric matrix ${\mathbf A}$ is similar to a
nonnegative matrix, is quite analogical to the proof of theorem 3.
Let us prove the irreducibility. Suppose the opposite: let the
nonnegative matrix $\widetilde{{\mathbf A}}$ be reducible, and the
initial ${\mathcal J}$--sign-symmetric matrix ${\mathbf A}$ be
irreducible. As it follows from the definition of reducibility,
there exists such a permutation of coordinates, that the matrix
$\widetilde{{\mathbf A}}$ will be reduced to canonical form (3).
Since the matrices ${\mathbf A}$ and $\widetilde{{\mathbf A}}$ are
connected by the similarity transformation with a diagonal matrix
${\mathbf D}$, then the matrix ${\mathbf A}$ will also be reduced
to form (3) by the same permutation of coordinates. We came to the
contradiction, because of the matrix ${\mathbf A}$ is irreducible.
$\square$

{\bf Corollary.} {\it Let the matrix $\mathbf A$ of a linear
operator $A$ be ${\mathcal J}$--sign-symmetric and irreducible.
Then the spectral radius $\rho(A)
> 0$ is a simple positive eigenvalue of the operator $A$. If $h$ is a
number of the eigenvalues of the operator $A$, which are equal in
modulus to $\rho(A)$, then all of them are simple and they
coincide with the $h$-th roots of $(\rho(A))^h$. More than that,
the spectrum of the operator $A$ is invariant under rotations by
$\frac{2\pi}{h}$ about the origin.}

It's easy to see, that the number of all different types of
strictly ${\mathcal J}$--sign-symmetric $n \times n$ matrices is
equal to the number of all subsets of the set $\{1, \ \ldots, \
n\}$, divided by $2$, i.e. $2^{n - 1}$. In turn, the number of all
different types of strictly ${\mathcal J}$--sign-symmetric $C_n^2
\times C_n^2$ matrices is equal to the number of all subsets of
the set $\{1, \ \ldots, \ C_n^2\}$, divided by $2$, i.e. $2^{C_n^2
- 1}$.

\section{The connection between the $W$--matrix and the second
compound matrix} Later on in this paper we shall study the case,
when the matrix ${\mathbf A}$ is ${\mathcal J}$--sign-symmetric,
i.e. is similar to a nonnegative matrix, and the second compound
matrix ${\mathbf A}^{(2)}$ is also ${\mathcal J}$--sign-symmetric,
i.e. is also similar to a non-negative matrix. Note, that this two
conditions absolutely do not mean, that the matrix ${\mathbf A}$
is similar to a $2$--totally positive matrix. This two conditions
also do not guarantee the reality of the peripheral spectrum of
the matrix ${\mathbf A}$. To show this, we should return to the
given above conception of a $W$--basis and a $W$--matrix. The
following theorem about the link between the structure of the
matrix ${\mathbf A}_W^{(2)}$, constructed with respect to an
arbitrary set $W$, which satisfies conditions (1) and (2), and the
structure of the second compound matrices ${\mathbf A}^{(2)}$, is
true.

{\bf Theorem 5.} {\it Let the second compound matrix ${\mathbf
A}^{(2)}$ of the matrix ${\mathbf A}$ be strictly ${\mathcal
J}$--sign-symmetric. Then there exists such a set $W \in \{1, \
\ldots, \ n\}\times \{1, \ \ldots, \ n\}$, which satisfies
conditions (1) and (2), that the corresponding $W$--matrix
${\mathbf A}_W^{(2)}$ is positive.

The inverse statement is also true. Let a $W$--matrix ${\mathbf
A}_W^{(2)}$ of a matrix ${\mathbf A}$ be positive. Then the second
compound matrix ${\mathbf A}^{(2)}$ is strictly ${\mathcal
J}$--symmetric.}

{\bf Proof.} $\Leftarrow$ Let a $W$--matrix ${\mathbf A}_W^{(2)}$
of a matrix ${\mathbf A}$, constructed with respect to a set $W$,
which satisfies conditions (1) and (2), be positive. Show, that
the second compound matrix ${\mathbf A}^{(2)}$ is strictly
${\mathcal J}$--sign-symmetric. Define the set ${\mathcal J}$ by
the following way:
$${\mathcal J} = \{\alpha: \ (i,j)\in (M \cap W) \setminus \Delta
\},$$ where $\alpha$ is the number of the set $(i,j)$ in the
lexicographic numeration. Verify, that the minor $A\begin{pmatrix}
  i & j \\
  k & l
\end{pmatrix}$, where $i < j$, $k < l$ is negative if and only if
one of the numbers of the sets $(i,j)$, $(k,l)$ belongs to the set
${\mathcal J}$, and the other belongs to the set $\{1, \ \ldots, \
C_n^2\}\setminus {\mathcal J}$. First examine the following
decomposition of the set $M$:
$$M = (M\cap W) \cup (M \cap \widetilde{W}).$$ Then the following
representation is true:
$$M \times M = ((M \cap W) \times (M \cap W)) \cup ((M \cap W)
\times (M \cap \widetilde{W})) \cup $$ $$ \cup ((M \cap
\widetilde{W})\times (M \cap W)) \cup ((M \cap
\widetilde{W})\times (M \cap \widetilde{W})).$$

 Analyze an arbitrary minor
$A\begin{pmatrix}
  i & j \\
  k & l
\end{pmatrix}$, where $i < j$, $k < l$. One of the following four cases takes place:
\begin{enumerate}
\item[\rm 1.] Both the numbers of the sets
$(i,j)$, $(k,l)$ belong to the set ${\mathcal J}$. In this case
both the pairs $(i,j)$ and $(k,l)$ belong to the set $M \cap W$.
The minor $A\begin{pmatrix}
  i & j \\
  k & l
\end{pmatrix}$ coincides with an element of a positive matrix ${\mathbf A}_W^{(2)}$,
as it follows, it is positive.
\item[\rm 2.] Both the numbers of the sets
$(i,j)$, $(k,l)$ belong to the set $\{1, \ \ldots, \
C_n^2\}\setminus {\mathcal J}$. It's easy to see, that the set
$\{1, \ \ldots, \ C_n^2\}\setminus {\mathcal J}$ consists of the
numbers of those and only those pairs, which belong to the set $M
\cap \widetilde{W}$. As it follows, both the pairs $(i,j)$,
$(k,l)$ belong to the set $M \cap \widetilde{W}$. Let us change
places of the rows and the columns of the minor (the minor will
change its sign twice). We'll receive the equality
$A\begin{pmatrix}
  i & j \\
  k & l
\end{pmatrix} = A\begin{pmatrix}
  j & i \\
  l & k
\end{pmatrix}$. Since both the pairs $(j,i)$ and $(l,k)$ belong
to the set $M \cap W$, then the minor $A\begin{pmatrix}
  j & i \\
  l & k
\end{pmatrix}$ coincides with an element of the matrix ${\mathbf
A}_W^{(2)}$ and, as it's follows, it's also positive.

\item[\rm 3.] The number of the pair $(i,j)$ belongs to the set
${\mathcal J}$, and the number of the pair $(k,l)$ belongs to the
set $\{1, \ \ldots, \ C_n^2\}\setminus {\mathcal J}$. In this case
the pair $(i,j) \in M \cap W$, and the pair $(k,l) \in M \cap
\widetilde{W}$. As it follows from the equality $A\begin{pmatrix}
  i & j \\
  k & l
\end{pmatrix} = - A\begin{pmatrix}
  i & j \\
  l & k
\end{pmatrix}$, the minor $A\begin{pmatrix}
  i & j \\
  k & l
\end{pmatrix}$ is opposite to an element of a positive matrix ${\mathbf
A}_W^{(2)}$, that's why it's negative.

\item[\rm 4.] The case, when the pair $(i,j) \in M \cap \widetilde{W}$,
and the pair $(k,l) \in M \cap W$, is studied by analogy.
\end{enumerate}

$\Rightarrow$ Let us prove the inverse. Let the second compound
matrix ${\mathbf A}^{(2)}$ be strictly ${\mathcal
J}$--sign-symmetric. As it follows, there exists a subset
${\mathcal J} \subseteq \{1, \ \ldots, \ C_n^2\}$, for which
${\mathbf a}^{(2)}_{\alpha\beta} < 0$ if and only if one of the
numbers $\alpha$, $\beta$ belongs to the set ${\mathcal J}$, and
the other belongs to the set $\{1, \ \ldots, \ C_n^2\}\setminus
{\mathcal J}$. Define a set $W$ for which the corresponding
$W$--matrix ${\mathbf A}_W^{(2)}$ is positive, by the following
way: $(i,j) \in W$ if and only if one of the following two cases
takes place:
\begin{enumerate}
\item[\rm (a)] $i < j$, and the number $\alpha$ of the pair $(i,j)$
(in the lexicographic numeration), belongs to the set ${\mathcal
J}$;

\item[\rm (b)] $i > j$, and the number $\widetilde {\alpha}$
of the pair $(j,i)$ belongs to the set $\{1, \ \ldots, \
C_n^2\}\setminus {\mathcal J}$.\end{enumerate}
 It's easy to see, that the set $W$ satisfies conditions (1) and (2).
 The positivity of the corresponding $W$--matrix is proved
by analogy with the first part of the proof of the theorem.
$\square$

 Later on we shall impose on the matrix
${\mathbf A}^{(2)}$ the condition of ${\mathcal J}$--symmertricity
together with the condition of irreducibility.

{\bf Theorem 6.} {\it Let the second compound matrix ${\mathbf
A}^{(2)}$ of a matrix ${\mathbf A}$ be ${\mathcal
J}$--sign-symmetric. Then there exists such a set $W \in \{1, \
\ldots, \ n\}\times \{1, \ \ldots, \ n\}$, which satisfies
conditions (1) and (2), that the corresponding $W$--matrix
${\mathbf A}_W^{(2)}$ is nonnegative. More than that, if ${\mathbf
A}^{(2)}$ is irreducible, then ${\mathbf A}_W^{(2)}$ is also
irreducible.

The inverse is also true. Let a $W$--matrix ${\mathbf A}_W^{(2)}$
of the matrix ${\mathbf A}$ be nonnegative. Then the second
compound matrix ${\mathbf A}^{(2)}$ is ${\mathcal J}$--sign
symmetric, and if ${\mathbf A}_W^{(2)}$ is irreducible, then
${\mathbf A}^{(2)}$ is also irreducible}.

{\bf Proof}. $\Rightarrow$ Let us define the set $W$,
corresponding to the set ${\mathcal J}$, as it was shown above in
the proof of theorem 5: $(i,j) \in W$ if and only if one of the
following two cases takes place:
\begin{enumerate}
\item[\rm (a)] $i < j$ and the number $\alpha$ of the pair $(i,j)$
in the lexicographic numeration belongs to the set ${\mathcal J}$;

\item[\rm (b)] $i > j$ and the number $\widetilde {\alpha}$
of the "inverse" \ pair $(j,i)$ in the lexicographic numeration
belongs to the set $\{1, \ \ldots, \ C_n^2\}\setminus {\mathcal
J}$.\end{enumerate}

 The nonnegativity of the matrix
${\mathbf A}_W^{(2)}$ is proved by analogy with the proof of the
positivity of the matrix ${\mathbf A}_W^{(2)}$ in theorem 5. The
irreducibility of the matrix ${\mathbf A}_W^{(2)}$ is proved by
analogy with the proof of the irreducibility of the matrix
$\widetilde{{\mathbf A}}$ in theorem 4. $\square$

\section{ Generalization of the Gantmacher--Krein theorems to the
case of a positive matrix with a strictly ${\mathcal
J}$--sign-symmetric second compound matrix.} Let us prove the
following theorem about the spectral properties of a positive
matrix with a strictly ${\mathcal J}$--sign-symmetric second
compound matrix.

 {\bf Theorem 7.} {\it Let the matrix ${\mathbf A}$ of a linear operator
$A$ be positive. Let its second compound matrix ${\mathbf
A}^{(2)}$ be strictly ${\mathcal J}$--sign-symmetric. Then the
operator $A$ has the first positive simple eigenvalue $\lambda_1 =
\rho(A)$, and the second positive simple eigenvalue $\lambda_2$}:
$$ \lambda_1 > \lambda_2 > |\lambda_3| > \ldots > 0.$$

{\bf Proof.} Enumerate the eigenvalues of the matrix ${\mathbf A}$
in order of decrease of their modules (taking into account their
multiplicities):
$$|\lambda_{1}| \geq | \lambda_{2}| \geq |\lambda_{3}| \geq \ldots
\geq |\lambda_{n}|$$
 Applying the Perron theorem to the matrix ${\mathbf A}$,
we get: $\lambda_{1} = \rho(A)>0$ is a simple positive eigenvalue
of ${\mathbf A}$. Examine the second compound matrix ${\mathbf
A}^{(2)}$, which is strictly ${\mathcal J}$--sign-symmetric.
Applying theorem 5 to the matrix ${\mathbf A}^{(2)}$, we get, that
there exists such a set $W$, for which the $W$--matrix ${\mathbf
A}^{(2)}_W$ is positive. Applying the Perron theorem to the matrix
${\mathbf A}^{(2)}_W$, we get: $\rho({\mathbf A}^{(2)}_W) > 0$ is
a simple positive eigenvalue of ${\mathbf A}^{(2)}_W$.

 As it follows from the statement of theorem 2, a $W$--matrix ${\mathbf
A}^{(2)}_W$ has no other eigenvalues, except all the possible
products of the form $\lambda_{i}\lambda_{j}$, where $ i < j$.
Therefore $\rho({\mathbf A}^{(2)}_W)>0$ can be represented in the
form of the product $\lambda_{i}\lambda_{j}$ with some values of
the indices $i,j$, \ $i < j$. It follows from the facts that the
eigenvalues are numbered in a decreasing order, and there is only
one eigenvalue on the spectral circle $|\lambda| = \rho({\mathbf
A})$, that
 $\rho({\mathbf
A}^{(2)}_W) =  \lambda_{1}\lambda_{2}$.
  Therefore $\lambda_{2} = \frac{\rho({\mathbf
A}^{(2)}_W)}{\lambda_{1}}>0$. $\square$

Let us give some examples, which illustrate theorem 7.

 {\bf Example 3.} Let $${\mathbf A} = \begin{pmatrix}
  30 & 41 & 3 & 16 \\
  41 & 61 & 3 & 20 \\
  3 & 3 & 1 & 2 \\
  16 & 20 & 2 & 10
\end{pmatrix}$$
 Then the second compound matrix is the following:
 $${\mathbf A}^{(2)} = \begin{pmatrix}
   149 & -33 & -56 & -60 & -156 & 12 \\
   -33 & 21 & 12 & 32 & 34 & -10 \\
   -56 & 12 & 44 & 22 & 90 & -2 \\
   -60 & 32 & 22 & 52 & 62 & -14 \\
   -156 & 34 & 90 & 62 & 210 & -10 \\
   12 & -10 & -2 & -14 & -10 & 6 \
 \end{pmatrix}$$
 In this case the set ${\mathcal J}$ is equal to $\{1, \ 6\}$ or
$ \{2, \ 3, \ 4, \ 5\}$. According to theorem 7, the operator $A$
has the first and the second simple positive eigenvalues (they are
equal to 97,0688 and 4,16138 respectively).

{\bf Example 4.} Let $$A = \begin{pmatrix}
  2 & 5 & 4 & 3 \\
  3 & 36 & 25 & 12 \\
 3 & 25 & 18 & 9 \\
 3 & 12 & 9 & 6
\end{pmatrix}$$

 Then the second compound matrix is the following:

 $${\mathbf A}^{(2)} = \begin{pmatrix}
   57 & 38 & 15 & -19 & -48 & -27 \\
   35 & 24 & 9 & -10 & -30 & -18 \\
   9 & 6 & 3 & -3 & -6 & -3 \\
   -33 & -21 & -9 & 23 & 24 & 9 \\
   -72 & -48 & -18 & 24 & 72 & 42 \\
   -39 & -27 & -9 & 9 & 42 & 27 \
 \end{pmatrix}$$

The set ${\mathcal J}$ is equal to $\{1, \ 2, \ 3\}$ or $\{4, \ 5,
\ 6\}$.

According to theorem 7, the operator $A$ has the first positive
simple eigenvalue (which is equal to 58.5009) and the second
positive simple eigenvalue (which is equal to 3.09002).

The following statement can be easily proved by analogy with
theorem 7, using the corollary from theorem 3 instead of Perron's
theorem: {\it let the matrix ${\mathbf A}$ of a linear operator
$A$ be strictly ${\mathcal J}$--sign-symmetric together with its
second compound matrix ${\mathbf A}^{(2)}$. Then the operator $A$
has the first positive simple eigenvalue $\lambda_1 = \rho(A)$,
and the second positive simple eigenvalue $\lambda_2$}:
$$ \lambda_1 > \lambda_2 > |\lambda_3| > \ldots > 0.$$

\section{ Permutations and isomorphisms of the space ${\mathbb X}$}
It is well-known (see [1], p. 317, theorem 13), that if the matrix
$\mathbf{A}$ of a linear operator $A$ is nonnegative together with
its second compound matrix ${\mathbf{A}}^{(2)}$, then the two
largest in modulus eigenvalues of the operator $A$ are real and
nonnegative. However, if we change the nonnegativity of the second
compound matrix ${\mathbf{A}}^{(2)}$ by the nonnegativity of a
$W$-matrix ${\mathbf A}_W^{(2)}$, then the peripheral spectrum of
the operator $A$ will not be always real. Later on we'll show,
that the reality of the peripheral spectrum of $A$ depends on if
the binary relation, defined by the set $W$, possesses the
additional property of transitivity.

Let us call the set $W$, which satisfies conditions (1) and (2),
{\it transitive}, if the binary relation, defined by $W$ on the
set $\{1, \ \ldots, \ n\}$, is a linear order relation. To analyze
the case, when the set $W$ is transitive, we shall use the
following lemma.

{\bf Lemma 2.} {\it Every set $W$, which defines a linear order
relation on the set $\{1, \ \ldots, \ n\}$, is uniquely defined by
the permutation $\theta = (\theta(1), \ldots, \theta(n))$. The
inverse is also true: every permutation $\theta$ of indices $\{1,
\ \ldots, \ n\}$ is uniquely defined by the set $W$, which defines
a linear order relation on the set $\{1, \ \ldots, \ n\}$}.

{\bf Proof}. $\Rightarrow$ Given a permutation $\theta =
(\theta(1), \ldots, \theta(n))$. Define a set $W$ by the following
way: the pair $(i,j) \in W$, if and only if $\theta^{-1}(i) \leq
\theta^{-1}(j)$, where $\theta^{-1}$ is the inverse permutation to
$\theta$. Show, that the set $W$ defines a linear order relation
on $\{1, \ \ldots, \ n\}$. Properties (1) and (2) are obvious. Let
us verify the property of transitivity. Let the inclusions $(i,j)
\in W$ and $(j,k) \in W$ be true for some indices $i, j, k \in
\{1, \ \ldots, \ n\}$. Then the inequalities $\theta^{-1}(i) \leq
\theta^{-1}(j)$ and $\theta^{-1}(j) \leq \theta^{-1}(k)$ are true.
As it follows from this two inequalities and the property of
transitivity of the natural linear order relation on $\{1, \
\ldots, \ n\}$, the inequality $\theta^{-1}(i) \leq
\theta^{-1}(k)$ and the inclusion $(i,k) \in W$ are true. The
transitivity is also realized.

 $\Leftarrow$ Given a set $W$, which defines a linear order relation
on $\{1, \ \ldots, \ n\}$. Let us define the permutation $\theta$
with the help of the following algorithm:

1) On the first step we define $\theta_1(1) = 1$.

2) On the second step we define $\theta_2(1) = 2$, $\theta_2(2) =
1$, if $(2,1) \in W $ and $\theta_2(1) = 1$, $\theta_2(2) = 2$
otherwise.

3) On the $j$th step we have a permutation $\theta_{j-1}$ of $j-1$
indices. Define $
 l = \max\{k: \ 1 \leq k \leq j-1; \
(\theta_{j - 1}(k), j) \in W \}$ (in the case, when for any $k: \
1 \leq k \leq j-1$ the inclusion $(\theta_{j - 1}(k), j) \in
\widetilde{W}$ is true, we define $l = 0$) and let
$$ \theta_j(i)
 = \left\{\begin{array}{cc} \theta_{j - 1}(i), &
\mbox{if $i \leq l$ };
\\[10pt] j , & \mbox{if $i = l + 1$ };
\\[10pt]\theta_{j - 1}(i - 1), & \mbox{if $l + 1 \leq i \leq j$ .}\end{array}\right.$$

On the $n$-th step we have a permutation of $n$ indices. Show,
that such a permutation defines the given set $W$. Define the set
$V$ by this permutation, as it was shown above in the first part
of the proof. Show, that the set $V$ coincides with the set $W$.
Let $(i,j) \in V$. In this case it follows from the inequality
$\theta^{-1}(i) \leq \theta^{-1}(j)$, that the index $i$ precedes
the index $j$ in the image of the set $\{1, \ \ldots, \ n \}$ by
the permutation $\theta$. Let $k_1, \ \ldots, \ k_m$ be indices,
disposed between $i$ and $j$ in $\theta(1, \ \ldots, \ n )$. Write
$\theta(1, \ \ldots, \ n )$ in the following form:
$$\theta(1, \ \ldots, \ n ) = \ldots, \ i, \ k_1, \ \ldots, \ k_m,
\ j, \ \ldots. $$ It follows from the construction of the
permutation $\theta$, that all the pairs $(i, k_1)$, $(k_2, k_3),
\ \ldots,$ $(k_{m-1}, k_m)$, $(k_m, j)$ belong to the set $W$.
Since $W$ is transitive, the inclusion $(i, k_2) \in W$ follows
from the inclusions $(i, k_1) \in W, \ (k_1, k_2) \in W$. Then,
the inclusion $(i, k_3) \in W$ follows from the inclusions $(i,
k_2) \in W, \ (k_2, k_3) \in W$. Repeating this reasoning for $m$
times, we'll get the inclusion $(i,j) \in W$. As it follows, the
inclusion $V \subseteq W$ is true. Prove the inverse inclusion.
Let $(i,j) \in W$. Prove, that $(i,j) \in V$. It's enough for this
to show, that $\theta^{-1}(i) \leq \theta^{-1}(j)$. Suppose the
opposite: let $\theta^{-1}(i)
> \theta^{-1}(j)$. Then the index $j$ precedes the index $i$
in the image of the set $\{1, \ \ldots, \ n \}$ by the permutation
$\theta$, and it follows from the above reasoning, that $(j,i) \in
W$. As it follows, both the pairs $(i,j)$ and $(j,i)$ belong to
the set $W$. This contradicts condition (2). $\square$

Let $Q_\theta$ be a permutation operator, defined on the basic
vectors by the following way: \linebreak $Q_\theta(e_i) =
e_{\theta(i)} \ \ (i = 1, \ \ldots, \ n)$. The following theorem
is true.

{\bf Theorem 8.} {\it Let the matrix ${\mathbf A}$ of a linear
operator $A: {\mathbb R}^n \rightarrow {\mathbb R}^n$ be
nonnegative, and let its second compound matrix ${\mathbf
A}^{(2)}$ be ${\mathcal J}$--sign-symmetric. Let the set $W$,
corresponding to the set of the indices ${\mathcal J}$, be
transitive. Then there exists such a permutation operator
$Q_\theta$, that the matrix ${\mathbf P} = {\mathbf
Q}^{T}_\theta{\mathbf A}{\mathbf Q}_\theta$ is nonnegative
together with its second compound matrix ${\mathbf P}^{(2)}$. More
than that, if the matrices ${\mathbf A}$ and ${\mathbf A}^{(2)}$
are irreducible, the matrices ${\mathbf P}$ and ${\mathbf
P}^{(2)}$ are also irreducible.}

{\bf Proof}. Define the permutation $\theta$ with respect to the
set $W$, using the algorithm, given above. It's easy to see, that
$p_{ij} = a_{\theta(i)\theta(j)}$. The matrix ${\mathbf P}
={\mathbf Q}^T_\theta{\mathbf A}{\mathbf Q}_\theta$ is obviously
nonnegative and irreducible. Prove the nonnegativity of the second
compound matrix ${\mathbf P}^{(2)}$. Study an arbitrary minor
$P\begin{pmatrix}
  i & j \\
  k & l
\end{pmatrix}$, where $i < j$, $k < l$. It is equal to
the generalized minor $A\begin{pmatrix}
  \theta(i) & \theta(j) \\
 \theta(k) & \theta(l)
\end{pmatrix}$.

It follows from the inequalities $i < j$, $k < l$ and the
construction of the permutation $\theta$, that both the pairs
$(\theta(i), \theta(j))$ and $(\theta(k), \theta(l))$ belong to
the set $W$ (a pair $(\theta(i),\theta(j)) \in W$ if and only if
$\theta^{-1}\theta(i) \leq \theta^{-1}\theta(j)$). As it follows,
the minor $A\begin{pmatrix}
  \theta(i) & \theta(j) \\
 \theta(k) & \theta(l)
\end{pmatrix}$ is an element of the $W$--matrix ${\mathbf
A}_W^{(2)}$, defined by the matrix ${\mathbf A}$ and the set $W$.
So it's easy to see, that the matrix ${\mathbf P}^{(2)}$ coincide
(up to a permutation of coordinates) with the $W$--matrix
${\mathbf A}_W^{(2)}$, which is, according to theorem 6,
nonnegative and irreducible.
 $\square$

Note, that in the case, when the set $W$ is not transitive, the
statement of theorem 8 is not always correct.

\section{Approximation of a ${\mathcal J}$--sign-symmetric matrix
by strictly ${\mathcal J}$--sign-symmetric matrices.}  Let us
prove the generalization of the statement about the approximation
of a totally positive matrix by strictly totally positive matrices
(see [1], p. 317, corollary from the theorem 12), using theorem 8.

{\bf Theorem 9.} {\it Let ${\mathbf A}$ be a nonnegative matrix.
Let its second compound matrix ${\mathbf A}^{(2)}$ be ${\mathcal
J}$--sign-symmetric. Let the set $W$, corresponding to one of
possible sets of the indices ${\mathcal J}$, be transitive. Then
there exists a sequence $\{{\mathbf A}_n\}$ of positive matrices
with strictly ${\mathcal J}$--sign-symmetric second compound
matrices, which converges to ${\mathbf A}$}.

{\bf Proof.} According to theorem 8, there exists a permutation
operator $Q_\theta$, for which the matrix ${\mathbf P} = {\mathbf
Q}^{T}_\theta{\mathbf A}{\mathbf Q}_\theta$ is nonnegative
together with its second compound matrix ${\mathbf P}^{(2)}$.
Using the statement about the approximation of a totally positive
matrix by strictly totally positive matrices (see [1], p. 317,
corollary from the theorem 12), let us construct a sequence of
positive matrices ${\mathbf P}_n$ with positive second compound
matrices ${\mathbf P}_n^{(2)}$, which converges to ${\mathbf P}$.
Examine the sequence ${\mathbf A}_n = {\mathbf Q}_\theta{\mathbf
P}_n{\mathbf Q}^T_\theta$. It's easy to see, that the sequence
$\{{\mathbf A}_n\}$ converges to the matrix ${\mathbf A}$. It's
also easy to see, that every matrix ${\mathbf A}_n^{(2)}$ is
strictly ${\mathcal J}$--sign-symmetric. $\square$

We get the following statement, using theorem 9 and the property
of continuity.

{\bf Theorem 10.} {\it Let ${\mathbf A}$ be a nonnegative matrix.
Let its second compound matrix ${\mathbf A}^{(2)}$ be ${\mathcal
J}$--sign-symmetric. Let the set $W$, corresponding to one of
possible sets of the indices ${\mathcal J}$, be transitive. Then
the two largest in modulus eigenvalues of the operator $A$ are
nonnegative.}

{\bf Proof.} The proof of theorem 10 follows from the statement of
theorem 7 of the existence of the two largest in modulus positive
simple eigenvalues of a positive matrix with a strictly ${\mathcal
J}$--sign-symmetric second compound matrix. $\square$

Later on we'll show, that if the set $W$, corresponding to the set
of the indices ${\mathcal J}$, is not transitive, then the
statement of theorem 10 of reality of the first two eigenvalues
will not be true.

Let us generalize the statements, proved above, to the class of
${\mathcal J}$--sign-symmetric matrices with ${\mathcal
J}$--sign-symmetric second compound matrices. Let ${\mathbf A}$ be
a ${\mathcal J}$--sign-symmetric matrix, and let ${\mathcal J}$ be
a subset of the set $\{1, \ \ldots, \ n\}$ in the definition of
${\mathcal J}$--sign-symmetricity (i.e. such a subset, that the
inequality $a_{ij} \leq 0$ is true for any two numbers $i,j$, one
of which belongs to the set ${\mathcal J}$, and the other belongs
to the set $\{1, \ \ldots, \ n\}\setminus {\mathcal J}$; and the
strict inequality $a_{ij} < 0$ is true only if one of the numbers
$i$, $j$ belongs to  ${\mathcal J}$, and the other belongs to
$\{1, \ \ldots, \ n\}\setminus {\mathcal J}$). Let ${\mathbf
A}^{(2)}$ be a ${\mathcal J}$--sign-symmetric matrix. Let
$\widetilde{{\mathcal J}}$ be a subset of $\{1, \ \ldots, C_n^2\}$
in the definition of ${\mathcal J}$--sign-symmetricity for the
matrix ${\mathbf A}^{(2)}$. Let us construct a set
$\widehat{W}({\mathcal J}, \widetilde{{\mathcal J}}) \subseteq
(\{1, \ \ldots, \ n\} \times \{1, \ \ldots, \ n\})$ with respect
to the sets ${\mathcal J}$ and $\widetilde{{\mathcal J}}$ by the
following way.

A pair $(i,j)$ belongs to the set $\widehat{W}({\mathcal J},
\widetilde{{\mathcal J}})$ if and only if one of the following
four cases takes place:
\begin{enumerate}
\item[\rm (a)] $i < j$, both the numbers $i,j$ belong either to the set ${\mathcal
J}$, or to the set $\{1, \ \ldots, \ n\}\setminus {\mathcal J}$,
and the number $\alpha$, corresponding to the pair $(i,j)$ in the
lexicographic numeration, belongs to the set $\widetilde{{\mathcal
J}}$;

\item[\rm (b)] $i < j$, one of the numbers $i,j$ belongs to the set ${\mathcal
J}$, and the other belongs to the set $\{1, \ \ldots, \
n\}\setminus {\mathcal J}$, and the number $\alpha$, corresponding
to the pair $(i,j)$ in the lexicographic numeration, belongs to
the set $\{1, \ \ldots, \ C_n^2\}\setminus \widetilde{{\mathcal
J}}$;

\item[\rm (c)] $i > j$, both the numbers $i,j$ belong either to the set ${\mathcal
J}$, or to the set $\{1, \ \ldots, \ n\}\setminus {\mathcal J}$,
and the number $\alpha$, corresponding to the pair $(j, i)$ in the
lexicographic numeration, belongs to the set $\{1, \ \ldots, \
C_n^2\}\setminus \widetilde{{\mathcal J}}$;

\item[\rm (d)]  $i > j$, one of the numbers $i,j$ belongs to the set ${\mathcal
J}$, the other belongs to the set $\{1, \ \ldots, \ n\}\setminus
{\mathcal J}$, and the number $\alpha$, corresponding to the pair
$(j, i)$ in the lexicographic numeration, belongs to the set
$\widetilde{{\mathcal J}}$.

\end{enumerate}

Let us prove the following statement.

{\bf Theorem 11.} {\it Let ${\mathbf A}$ be a ${\mathcal
J}$--sign-symmetric matrix. Let its second compound matrix
${\mathbf A}^{(2)}$ also be ${\mathcal J}$--sign-symmetric. Let
the set $\widehat{W}({\mathcal J}, \widetilde{{\mathcal J}})$ be
transitive. Then there exists a sequence $\{{\mathbf A}_n\}$ of
strictly ${\mathcal J}$--sign-symmetric matrices with strictly
${\mathcal J}$--sign-symmetric second compound matrices, which
converges to ${\mathbf A}$.}

{\bf Proof.} Let ${\mathbf A}$ be a ${\mathcal J}$--sign-symmetric
matrix. Then, according to theorem 4, the matrix ${\mathbf A}$ can
be represented in the form:
$${\mathbf A} = {\mathbf D} \widetilde{{\mathbf A}} {\mathbf
D}^{-1}, \eqno(4)$$ where $\widetilde{{\mathbf A}}$ is a
nonnegative matrix. Examine the second compound matrix ${\mathbf
A}^{(2)}$. It's known (see [1], p. 80, property 1), that the
compound matrix of a product of matrices, is equal to the product
of compound matrices of the factors. As it follows, the matrix
${\mathbf A}^{(2)}$ can be represented in the form:
$${\mathbf A}^{(2)} = {\mathbf D}^{(2)} \widetilde{{\mathbf A}}^{(2)}
 ({\mathbf D}^{-1})^{(2)}. $$
It's also known (see [1], p. 80, property 2), that the inverse
matrix of the compound matrix is equal to the compound matrix of
the inverse matrix. I.e. $({\mathbf D}^{-1})^{(2)} = ({\mathbf
D}^{(2)})^{-1}$, and the equality
$${\mathbf A}^{(2)} = {\mathbf D}^{(2)} \widetilde{{\mathbf A}}^{(2)}
 ({\mathbf D}^{(2)})^{-1} \eqno(5)$$ is correct.

Express the matrix $\widetilde{{\mathbf A}}^{(2)}$ from equality
(5):
$$\widetilde{{\mathbf A}}^{(2)} = ({\mathbf D}^{(2)})^{-1} {\mathbf
A}^{(2)}
 {\mathbf D}^{(2)}. \eqno(6)$$ Both the matrices ${\mathbf D}^{(2)}$
and $({\mathbf D}^{(2)})^{-1}$ are diagonal matrices, which
diagonal elements are equal to $\pm 1$. According to the
conditions of the theorem, the matrix ${\mathbf A}^{(2)}$ is
${\mathcal
 J}$--sign-symmetric. So, it's easy to see,
that the matrix $\widetilde{{\mathbf A}}^{(2)}$ is also ${\mathcal
J}$--sign-symmetric.
 As it follows, we can apply given above
theorem 9 to the nonnegative matrix $\widetilde{{\mathbf A}}$ with
a ${\mathcal J}$--sign-symmetric second compound matrix
$\widetilde{{\mathbf A}}^{(2)}$. According to theorem 9, if the
set $W$, corresponding to the set $\widehat{{\mathcal J}}$ in the
definition of ${\mathcal J}$--sign-symmetricity of the matrix
$\widetilde{{\mathbf A}}^{(2)}$, is transitive, then there exists
a sequence $\{\widetilde{{\mathbf A}}_n\}$ of positive matrices
with strictly ${\mathcal J}$--sign-symmetric second compound
matrices, which converges to the matrix $\widetilde{{\mathbf A}}$.
Let us construct the sequence $\{{\mathbf A}_n\}$ by the following
way: ${\mathbf A}_n = {\mathbf D} \widetilde{{\mathbf A}}_n
{\mathbf D}^{-1},$ where ${\mathbf D}$ is a diagonal matrix in
equality (4). It's easy to see, that the sequence $\{{\mathbf
A}_n\}$ converges to the matrix ${\mathbf A}$. It's also easy to
see, that for any $n = 1, 2, \ldots$ both the matrix ${\mathbf
A}_n$ and its second compound matrix are strictly ${\mathcal
J}$--sign-symmetric.

 Show, that the set $W$,
corresponding to the set $\widehat{{\mathcal J}}$ in the
definition of ${\mathcal J}$--sign-symmetricity of the matrix
$\widetilde{{\mathbf A}}^{(2)}$ coincide with
$\widehat{W}({\mathcal J}, \widetilde{{\mathcal J}})$. Write the
matrix $\widetilde{{\mathbf A}}^{(2)}$, using equality (6):
$$\widetilde{{\mathbf A}}^{(2)} = ({\mathbf D}^{(2)})^{-1} {\mathbf
A}^{(2)}
 {\mathbf D}^{(2)}. $$ Applying theorem 4 to the matrix
${\mathbf A}^{(2)}$, we'll get: $${\mathbf A}^{(2)} =
\widehat{{\mathbf D}}\widehat{{\mathbf A}} \widehat{{\mathbf
D}}^{-1},$$ where $\widehat{{\mathbf A}}$ is a nonnegative
 $C_n^2 \times C_n^2$ matrix, $\widehat{{\mathbf D}}$
is a diagonal matrix, which diagonal elements are equal to $\pm
1$. Therefore, the following equality is true:
$$\widetilde{{\mathbf A}}^{(2)} = ({\mathbf D}^{(2)})^{-1}
\widehat{{\mathbf D}}\widehat{{\mathbf A}} \widehat{{\mathbf
D}}^{-1} {\mathbf D}^{(2)}. \eqno (7)$$

Write equality (7) in the following form:
$$\widetilde{{\mathbf A}}^{(2)} =
\widetilde{{\mathbf D}}\widehat{{\mathbf A}} \widetilde{{\mathbf
D}}^{-1},$$ where $\widetilde{{\mathbf D}} = ({\mathbf
D}^{(2)})^{-1} \widehat{{\mathbf D}}$. Because of ${\mathbf
D}^{(2)}$ is a diagonal matrix, which diagonal elements are equal
to $\pm 1$, it's obvious, that $({\mathbf D}^{(2)})^{-1} =
{\mathbf D}^{(2)}$, and, as it follows, $\widetilde{{\mathbf D}} =
{\mathbf D}^{(2)} \widehat{{\mathbf D}}$.

It follows from the proofs of theorems 3 and 4, that the set
$\widehat{{\mathcal J}}$ is defined by the matrix
$\widetilde{{\mathbf D}}$ by the following way: the index $\alpha$
belongs to the set $\widehat{{\mathcal J}}$ if and only if the
element $\widetilde{d}_{\alpha\alpha} = - 1$. It's also easy to
see, that the elements of the diagonal matrix ${\mathbf D}^{(2)}$
can be defined by the set ${\mathcal J}$ by the following way:

$d^{(2)}_{\alpha\alpha}
 =  - 1$, if one of the indices $(i,j)$ of the pair
with the number $\alpha$ in the lexicographic numeration, belongs
to the set ${\mathcal J}$, and the other belongs to $\{1, \
\ldots, \ n\}\setminus {\mathcal J}$;

$d^{(2)}_{\alpha\alpha}
 =   1$, if both the indices $(i,j)$ belong either to the set
${\mathcal J}$ or to the set $\{1, \ \ldots, \ n\}\setminus
{\mathcal J}$.

The elements of the diagonal matrix $\widehat{{\mathbf D}}$ are
defined by the following way:
$$\widehat{d}_{\alpha\alpha}
 = \left\{\begin{array}{cc} - 1, &
\mbox{if $\alpha \in \widetilde{{\mathcal J}}$ };
\\[10pt] 1, &  \mbox{otherwise.}\end{array}\right.$$

The equality $\widetilde{d}_{\alpha\alpha} =
d^{(2)}_{\alpha\alpha}\widehat{d}_{\alpha\alpha}$ is true for the
elements of the matrix $\widetilde{{\mathbf D}}$. As it follows,
the elements of the set $\widehat{{\mathcal J}}$ can be defined by
the following way. The element $\alpha$ belongs to the set
$\widehat{{\mathcal J}}$ if and only if one of the following two
cases takes place:
\begin{enumerate}
\item[\rm (a)]both the numbers $i,j$ of the pair with the number
$\alpha$ in the lexicographic numeration, belongs either to the
set ${\mathcal J}$ or to the set $\{1, \ \ldots, \ n\}\setminus
{\mathcal J}$, and the number $\alpha$ belongs to the set
$\widetilde{{\mathcal J}}$;

\item[\rm (b)] one of the numbers $i,j$ belongs to the set ${\mathcal
J}$, the other belongs to the set $\{1, \ \ldots, \ n\}\setminus
{\mathcal J}$, and the number $\alpha$ belongs to the set $\{1, \
\ldots, \ C_n^2\}\setminus \widetilde{{\mathcal J}}$.
\end{enumerate}

It's easy to see, that the set $W$, corresponding to such a set
$\widehat{{\mathcal J}}$, will coincide with
$\widehat{W}({\mathcal J}, \widetilde{{\mathcal J}})$. $\square$

The following statement comes out from theorem 11.

{\bf Theorem 12.} {\it Let ${\mathbf A}$ be a ${\mathcal
J}$--sign-symmetric matrix. Let its second compound matrix
${\mathbf A}^{(2)}$ be also ${\mathcal J}$--sign-symmetric. Let
the set $\widehat{W}({\mathcal J}, \widetilde{{\mathcal J}})$ be
transitive. Then the two largest in modulus eigenvalues of the
operator $A$ are nonnegative.}

Note, that if the set $\widehat{W}({\mathcal J},
\widetilde{{\mathcal J}})$ is not transitive, then the
approximation of a ${\mathcal J}$--sign-symmetric matrix with a
${\mathcal J}$--sign-symmetric second compound matrix by strictly
${\mathcal J}$--sign-symmetric matrices with strictly ${\mathcal
J}$--sign-symmetric second compound matrices is not always
possible, and the statement of the reality of the two largest in
modulus eigenvalues is not always true.

\section{ Generalization of the Gantmacher--Krein theorems to the
case of an irreducible nonnegative matrices with an irreducible
${\mathcal J}$--sign-symmetric second compound matrix.}
 One can receive more detailed information on
the structure of the spectrum, than given in theorems 10 and 12,
imposing on the matrix ${\mathbf A}$ and it second compound matrix
${\mathbf A}^{(2)}$ the additional condition of irreducibility.

{\bf Theorem 13.} {\it Let the matrix ${\mathbf A}$ of a linear
operator $A: {\mathbb R}^n \rightarrow {\mathbb R}^n$ be
nonnegative and irreducible and let its second compound matrix
${\mathbf A}^{(2)}$ be irreducible and ${\mathcal
J}$--sign-symmetric. Let the set $W$, corresponding to the set of
the indices ${\mathcal J}$, be transitive. Then $h(A) = 1$ and the
operator $A$ has the first simple positive eigenvalue $\lambda_1$,
equal to the spectral radius $\rho(A)$, and the second simple
positive eigenvalue $\lambda_2$:
$$ \lambda_1
> \lambda_2 \geq |\lambda_3|
> \ldots.$$  If $h(A) = h(A \wedge A) = 1$, then $\lambda_{2}$
is different in modulus from the other eigenvalues. If $h(A) = 1$,
and $h(A \wedge A) > 1$, then the operator $A$ has $h(A \wedge A)$
eigenvalues $\lambda_{2}, \lambda_{3}, \ldots, \lambda_{h(A \wedge
A) + 1} $, equal in modulus to $\lambda_{2}$, each of them is
simple, and they coincide with the $h(A \wedge A)$th roots of
$\lambda_{2}^{h(A \wedge A)}$.}

{\bf Proof.} The equality $h(A) = 1$ follows from the statement of
theorem 10 of the reality of the peripheral spectrum. The
existence and the positivity of the first and the second
eigenvalues are proved by analogy with theorem 7. The simplicity
of $\lambda_{2}$ follows from the equality $\lambda_{2} =
\frac{\rho(A\wedge A)}{\rho(A)}$ and the simplicity of the
eigenvalues $\rho(A)$ and $\rho(A\wedge A)$.

 In the case of $h(A) = h(A \wedge A) = 1$
the distinction in modulus of $\lambda_{2}$ from the other
eigenvalues is obvious.

In the case of $h(A \wedge A) > 1$ it follows from theorem 2 and
the properties of the peripheral spectrum of the imprimitive
operator $A \wedge A$, that the equality $\lambda_{j}=
\frac{\rho(A \wedge A) e^\frac{2\pi(j-1)i}{h(A \wedge
A)}}{\rho(A)}$ is true for the eigenvalues $\lambda_{j}, \ j = 2,
\ \ldots, \ h(A \wedge A) + 1$. $\square$

Later on we shall show, that if the set $W$ is not transitive,
then the spectrum of the operator $A$ has another structue (there
is just three eigenvalues on the spectral circle $|\lambda| =
\rho(A)$).

{\bf Theorem 14.} {\it Let the matrix ${\mathbf A}$ of a linear
operator $A: {\mathbb R}^n \rightarrow {\mathbb R}^n$ be
nonnegative and irreducible, and let its second compound matrix
${\mathbf A}^{(2)}$ be ${\mathcal J}$--sign-symmetric and also
irreducible. Let the set $W$, defined by the set of the indices
${\mathcal J}$ be non-transitive. Then the operator $A$ has the
first positive eigenvalue $\lambda = \rho(A)$ with corresponding
positive eigenvector $x_1$. More than that, there is just three
eigenvalues on the spectral circle $|\lambda| = \rho(A)$, all of
them are simple and coincide with $3$th roots of $(\rho(A))^3$.
The following equality for the indices of imprimitivity of the
operators $A$ and $A \wedge A$ is true: $h(A) = h(A \wedge A) =
3$.}

{\bf Proof} Let us prove that $h(A) = h(A \wedge A) = 3$ by
contradiction, excluding all the possible values $h(A)$, except
$h(A) = 3$.

Suppose $h(A) = 1$. Enumerate the eigenvalues of the operator $A$,
repeated according to multiplicity, in order of decrease of their
modules:
$$|\lambda_{1}| \geq | \lambda_{2}| \geq  \ldots \geq
|\lambda_n|.$$ Applying the Frobenius theorem to the matrix
${\mathbf A}$, we'll get, that the operator $A$ has the first
positive eigenvalue $\lambda_{1} = \rho(A)>0$ with the
corresponding positive eigenvector $x_1$, and $\lambda_{1}$ is
simple and different in modulus from the other eigenvalues.
Applying the Frobenius theorem to the matrix ${\mathbf
A}_W^{(2)}$, which is also nonnegative and irreducible, we'll get,
that $\rho(A \wedge A)$ is a simple positive eigenvalue of the
operator $A \wedge A$, with the corresponding positive eigenvector
$\varphi$.

It follows from the fact, that there is only one eigenvalue on the
spectral circle $|\lambda| = \rho(A)$, and the statement of
theorem 2, that the eigenvalue $\rho(A \wedge A)
> 0$ can be represented in the form $\lambda_{1}\lambda_{m}$,
with some unique value $m > 1$. Without loss of generality of the
reasoning, we'll assume $m = 2$, i.e., that $\rho(A \wedge A)=
\lambda_{1}\lambda_{2}$. Then the eigenvector $\varphi$,
corresponding to the eigenvalue $\rho(A \wedge A)$ can be
represented in the form of the exterior product $x_1 \wedge x_2$
of the eigenvector $x_1$, corresponding to the eigenvalue
$\lambda_{1}$, and the eigenvector $x_2$, corresponding to the
eigenvalue $\lambda_{2}$. Let us examine the coordinates of the
vector $\varphi$ in the $W$--basis, defined by the set $W$. All of
them are positive. Because of $W$ is not transitive, there exists
at least one triple of indices $i,j,k \in \{1, \ \ldots, \ n \}$,
for which the inclusions $(i,j), \ (j,k) \in W$, $(i,k) \in
\widetilde{W}$ are true. In this case we'll receive the following
system of inequalities for the coordinates of the vector $\varphi
= x_1 \wedge x_2$ in the $W$--basis:
$$\varphi_\alpha = x_i^1x_j^2 - x_j^1x_i^2 > 0; $$
$$\varphi_\beta = x_j^1x_k^2 - x_k^1x_j^2 > 0; $$ $$\varphi_\gamma
= x_k^1x_i^2 - x_i^1x_k^2 > 0. $$ (Here $\alpha$, $\beta$,
$\gamma$ are the numbers of the pairs $(i,j)$, $(j,k)$ and $(k,i)$
respectively, $x_i^l$, $x_j^l$, $x_k^l$ are the coordinates of the
vectors $x_l$, $l = 1,2$). Let us multiply the first inequality by
$x_k^1$, and the second inequality by $x_i^1$. (Note, that all the
coordinates of the vector $x_1$ are positive, therefore the sings
of the inequalities will not change after multiplication.) Add
together both the inequalities. We'll get the system:
$$x_j^1(x_i^1x_k^2 - x_k^1x_i^2) > 0;
$$
$$ x_k^1x_i^2 - x_i^1x_k^2 > 0. $$
 It's easy to see, that this system has no solutions.
We came to the contradiction with the fact, that $\varphi = x_1
\wedge x_2$, where $x_1$ is a positive vector. so the case of
$h(A) = 1$ is excluded.

Exclude the case of $h(A) = 2$. For this we will prove that if a
nonnegative operator $A$ is imprimitive with $h(A) = 2$, then its
exterior square can not be nonnegative. Really, according to the
Frobenius theorem, there are two eigenvalues $\rho(A)
> 0$ and $- \rho(A)$ on the spectral circle $|\lambda| = \rho(A)$
of the operator $A$. As it follows, there is only one negative
eigenvalue $- \rho^{2}(A)$ on the spectral circle $|\lambda| =
\rho(A \wedge A)$ of the operator $A \wedge A$, and that is
impossible, if $A \wedge A$ is nonnegative.

Exclude the case of $h(A)>3$. Prove that the operator $A \wedge A$
is reducible, if it is nonnegative and $h(A)>3$. Really, it
follows from theorem 2 and the imprimitivity of $A$ that all the
eigenvalues of the operator $A \wedge A$ on the spectral circle
$|\lambda| = \rho(A \wedge A)$, can be represented as couple
products of different $h(A)$th roots of $\rho(A)^{h(A)}$. Let us
examine $\lambda_{j} = \rho(A)e^{\frac{2\pi(j-1)i}{h(A)}} \ (j =
1, \ldots, h(A))$ --- all the eigenvalues of the operator $A$,
situated on the spectral circle $|\lambda| = \rho(A)$. It's
obvious, that $\lambda_{2}\lambda_{h(A)} =
\lambda_{3}\lambda_{h(A)-1} = \ldots = \lambda_{k}\lambda_{h(A) -
(k - 2)} = \ldots = \rho(A)^{2}$. As it follows, the eigenvalue
$\rho(A \wedge A) = \rho(A)^{2}$ of the operator $A \wedge A$ is
not simple, and that is impossible, if $A \wedge A$ is
irreducible.

The only possible case is $h(A) = 3$. Let us prove, that if $A$ is
imprimitive with the index of imprimitivity $h(A) = 3$, and its
exterior square is irreducible, that $A \wedge A$ is also
imprimitive with $h(A \wedge A) = 3$. Really, in this case there
are three eigenvalues $\lambda_{1} = \rho(A)$, $\lambda_{2} =
\rho(A)e^{\frac{2\pi i}{3}}$, $\lambda_{3} = \rho(A)e^{\frac{4\pi
i}{3}}$ on the spectral circle $|\lambda| = \rho(A)$, and there is
also three eigenvalues $\lambda_{1}\lambda_{2} =
\rho(A)^{2}e^{\frac{2\pi i}{3}}$, $\lambda_{1}\lambda_{3} =
\rho(A)^{2}e^{\frac{4\pi i}{3}}$ and $\lambda_{2}\lambda_{3} =
\rho(A)e^{\frac{2\pi i}{3}} \rho(A)e^{\frac{4\pi i}{3}} =
\rho(A)^{2}$, which coincide with the 3th roots of
$(\rho(A)^{2})^{3}$, on the spectral circle $\lambda_{1} = \rho(A
\wedge A)$. $\square$

This statement follows from theorem 14: {\it if the matrix
${\mathbf A}$ of a linear operator $A: {\mathbb R}^n \rightarrow
{\mathbb R}^n$ is primitive, and its second compound matrix
${\mathbf A}^{(2)}$ is irreducible and ${\mathcal
J}$--sign-symmetric, then the set $W$, corresponding to the set of
the indices ${\mathcal J}$, is transitive.}

{\bf Example 5.} Let the operator $A:{\mathbb{R}}^{3} \rightarrow
{\mathbb{R}}^{3}$ be defined by the matrix $$\mathbf{A} =
\begin{pmatrix}
  0 & 0 & 1 \\
  1 & 0 & 0 \\
  0 & 1 & 0
\end{pmatrix}.$$
 This matrix is obviously nonnegative and irreducible.

In this case the second compound matrix is the following:
$${\mathbf A}^{(2)} =
\begin{pmatrix}
  0 & - 1 & 0 \\
  0 & 0 & - 1 \\
  1 & 0 & 0
\end{pmatrix}.$$

 The matrix ${\mathbf A}^{(2)}$ is obviously
${\mathcal J}$--sign-symmetric and irreducible. (In this case the
set ${\mathcal J}$ consists of two indices $1$ and $3$, or the one
index $2$). Examine the set $W$, corresponding to the set of
indices ${\mathcal J} = \{1, \ 3\}$. It consists of the pairs
$(1,2)$ and $(2,3)$, which have the numbers $1$ and $3$ in the
lexicographic numeration, and the pair $(3,1)$, the "inverse" \ of
which $(1,3)$ has the number $2$ (see illustration 1).

\begin{center}
\begin{picture}(70, 60)(125
, 25)

\put(125,25){\circle{5}}\put(125,50){\circle*{5}}\put(125,75){\circle{5}}
 \put(150,25){\circle{5}}\put(150,50){\circle{5}}\put(150,75){\circle*{5}}
 \put(175,25){\circle*{5}}\put(175,50){\circle{5}}\put(175,75){\circle{5}}
\end{picture}
\end{center}
\begin{center}
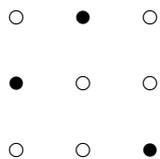
\small {Illustration 1. The set $W$.}
\end{center}
Such a set $W$ defines the non-transitive binary relation $(1
\prec 2), \ (2 \prec 3)$, $(3 \prec 1)$ on the set of the indices
$\{1, \ 2, \ 3 \}$. The operator $A$ satisfies the conditions of
theorem 14. It's easy to see, that $A$ has the first positive
simple eigenvalue $\lambda = \rho(A) = 1$, and there is just three
eigenvalues $1$, $e^{\frac{2\pi i}{3}}$ and $e^{\frac{4\pi i}{3}}$
on the spectral circle $|\lambda| = 1$, all of them are simple and
coincide with $3$th roots of $1$.

\section{ Generalization of the Gantmacher--Krein theorems to the
case of an irreducible ${\mathcal J}$--sign-symmetric matrix with
an irreducible ${\mathcal J}$--sign-symmetric second compound
matrix.}

Let us prove the generalization of the Gantmacher--Krein theorem
to the case of an irreducible ${\mathcal J}$--sign-symmetric
matrix, with an irreducible ${\mathcal J}$--sign-symmetric second
compound matrix.

{\bf Theorem 15.} {\it Let the matrix ${\mathbf A}$ of a linear
operator $A: {\mathbb R}^n \rightarrow {\mathbb R}^n$ be
${\mathcal J}$--sign-symmetric and irreducible. Let its second
compound matrix ${\mathbf A}^{(2)}$ be also ${\mathcal
J}$--sign-symmetric and irreducible. Then one of the following two
cases takes place:
\begin{enumerate}
\item[\rm (1)] The set $\widehat{W}({\mathcal J},
\widetilde{{\mathcal J}})$ is transitive. Then $h(A) = 1$, $h(A
\wedge A)$ is an arbitrary, and the operator $A$ has two positive
simple eigenvalues $\lambda_{1}$, $\lambda_{2}$:
$$\rho(A) = \lambda_{1} > \lambda_{2} \geq |\lambda_{3}| \geq
\ldots \geq |\lambda_{n}|. $$ If $h(A) = h(A \wedge A) = 1$, then
$\lambda_{2}$ is different in modulus from the other eigenvalues.
If $h(A) = 1$, and $h(A \wedge A)
> 1$, then the operator $A$ has $h(A \wedge A)$ eigenvalues $\lambda_{2}, \lambda_{3}, \ldots, \lambda_{h(A \wedge A) + 1}
$, equal in modulus to $\lambda_{2}$, each of them is simple, and
they coincide with the $h(A \wedge A)$th roots from
$\lambda_{2}^{h(A \wedge A)}$.
\item[\rm (2)] the set $\widehat{W}({\mathcal J},
\widetilde{{\mathcal J}})$ is not transitive. Then $h(A) = h(A
\wedge A) = 3$, and there is just three eigenvalues on the
spectral circle $|\lambda| = \rho(A)$. Each of them is simple, and
they coincide with the $3$th roots of $(\rho(A))^3$.
\end{enumerate}}

{\bf Proof.} Applying theorem 4, write the matrix ${\mathbf A}$ in
form (4):
$${\mathbf A} = {\mathbf D} \widetilde{{\mathbf A}} {\mathbf
D}^{-1}, $$ where $\widetilde{{\mathbf A}}$ is a nonnegative
irreducible matrix, ${\mathbf D}$ is a diagonal matrix, which
diagonal elements are equal to $\pm 1$. In this case the second
compound matrix ${\mathbf A}^{(2)}$ can be represented in form
(5):
$${\mathbf A}^{(2)} = {\mathbf D}^{(2)} \widetilde{{\mathbf A}}^{(2)}
 ({\mathbf D}^{(2)})^{-1}.$$
 It's
shown above in the proof of theorem 11, that the ${\mathcal
J}$--sign-symmetricity of the matrix $\widetilde{{\mathbf
A}}^{(2)}$ follows from equality (5).

According to the conditions of the theorem, the matrix ${\mathbf
A}^{(2)}$ is irreducible. So, it's easy to see, that the matrix
$\widetilde{{\mathbf A}}^{(2)}$ is also irreducible. As it
follows, we can apply given above theorems 13 and 14 to the
nonnegative irreducible matrix $\widetilde{{\mathbf A}}$ with a
${\mathcal J}$--sign-symmetric irreducible second compound matrix.
It follows from the similarity of the matrices ${\mathbf A}$ and
$\widetilde{{\mathbf A}}$, that their spectra coincide. I.e.
applying theorems 13 and 14 we'll get, that one of the following
two cases takes place:
\begin{enumerate}
\item[\rm (1)] The set $W$, corresponding to the set $\widehat{{\mathcal J}}$
in the definition of ${\mathcal J}$--sign-symmetricity of the
matrix $\widetilde{{\mathbf A}}^{(2)}$, is transitive. Then $h(A)
= 1$, $h(A \wedge A)$ is an arbitrary, and the operator $A$ has
two positive simple eigenvalues $\lambda_{1}$, $\lambda_{2}$:
$$\rho(A) = \lambda_{1}
> \lambda_{2} \geq |\lambda_{3}| \geq \ldots \geq |\lambda_{n}|.
$$ If $h(A) = h(A \wedge A) = 1$, then $\lambda_{2}$ is different
in modulus from the other eigenvalues. If $h(A) = 1$, and $h(A
\wedge A) > 1$, then the operator $A$ has $h(A \wedge A)$
eigenvalues $\lambda_{2}, \lambda_{3}, \ldots, \lambda_{h(A \wedge
A) + 1} $, equal in modulus to $\lambda_{2}$, each of them is
simple and they coincide with the $h(A \wedge A)$th roots of
$\lambda_{2}^{h(A \wedge A)}$.
\item[\rm (2)] The set $W$,
corresponding to the set of the indices $\widehat{{\mathcal J}}$,
is not transitive. Then $h(A) = h(A \wedge A) = 3$, and there is
just three eigenvalues on the spectral circle $|\lambda| =
\rho(A)$. All of them are simple and coincide with $3$th roots of
$(\rho(A))^3$.
\end{enumerate}

It's shown in the proof of theorem 11, that the set $W$,
corresponding to the set of the indices $\widehat{{\mathcal J}}$
in the definition of ${\mathcal J}$--sign-symmetricity of
$\widetilde{{\mathbf A}}^{(2)}$, coincide with the set
$\widehat{W}({\mathcal J}, \widetilde{{\mathcal J}})$. $\square$

{\bf Example 6.} Let the operator $A:{\mathbb{R}}^{3} \rightarrow
{\mathbb{R}}^{3}$ be defined by the matrix $$\mathbf{A} =
\begin{pmatrix}
  8,5 & 0 & 6,1 \\
  -5,6 & 3,2 & -7,4 \\
  6 & -2,8 & 6,6
\end{pmatrix}.$$
 This matrix is ${\mathcal J}$--sign-symmetric and irreducible.
In this case the set ${\mathcal J}$ in the definition of
${\mathcal J}$--sign-symmetricity of the matrix $\mathbf{A}$
consists of two indices $1$ and $3$.

In this case the second compound matrix is the following:
$${\mathbf A}^{(2)} =
\begin{pmatrix}
  27,2 & -28,74 & -19,52 \\
  -23,8 & 19,5 & 17,08 \\
  -3,52 & 7,44 & 0,4
\end{pmatrix}.$$

 The matrix ${\mathbf A}^{(2)}$ is also
${\mathcal J}$--sign-symmetric and irreducible. The set
$\widetilde{{\mathcal J}}$ consists of two indices $2$ and $3$.

 Examine the set $\widehat{W}({\mathcal J},
\widetilde{{\mathcal J}})$. The pair $(1,2)$ belongs to
$\widehat{W}({\mathcal J}, \widetilde{{\mathcal J}})$, because of
$1 < 2$, $1 \in {\mathcal J}$, $2 \in \{1, \ 2, \ 3\}\setminus
{\mathcal J}$, and the number $1$, corresponding to the pair
$(1,2)$ in the lexicographic numeration, belongs to the set $\{1,
\ 2, \ 3\}\setminus \widetilde{{\mathcal J}}$. The pair $(1,3)$
belongs to $\widehat{W}({\mathcal J}, \widetilde{{\mathcal J}})$,
because of $1 < 3$, both the numbers $1$ and $3$ belong to the set
${\mathcal J}$, and the number $2$, corresponding to the pair
$(1,3)$ in the lexicographic numeration, belongs to the set
$\widetilde{{\mathcal J}}$. And the pair $(3,2)$ belongs to
$\widehat{W}({\mathcal J}, \widetilde{{\mathcal J}})$, because of
$3 > 2$, $3 \in {\mathcal J}$, $2 \in \{1, \ 2, \ 3\}\setminus
{\mathcal J}$, and the number $3$, corresponding to the pair
$(2,3)$ in the lexicographic numeration, belongs to the set
$\widetilde{{\mathcal J}}$.

\begin{center}
\begin{picture}(70, 60)(125
, 25)

\put(125,25){\circle{5}}\put(125,50){\circle*{5}}\put(125,75){\circle*{5}}
 \put(150,25){\circle{5}}\put(150,50){\circle{5}}\put(150,75){\circle{5}}
 \put(175,25){\circle*{5}}\put(175,50){\circle{5}}\put(175,75){\circle{5}}
\end{picture}
\end{center}
\begin{center}
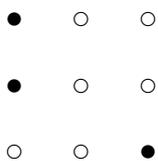
\small {Illustration 2. The set $\widehat{W}({\mathcal J},
\widetilde{{\mathcal J}})$.}
\end{center}
Such a set $\widehat{W}({\mathcal J}, \widetilde{{\mathcal J}})$
defines the linear order relation $1 \prec 3 \prec 2$ on the set
of the indices $\{1, \ 2, \ 3 \}$. The operator $A$ satisfies the
conditions of theorem 15, case (1). It's easy to see, that $A$ has
the first positive simple eigenvalue $\lambda_1 = \rho(A) =
15,102$, and the second positive simple eigenvalue $\lambda_{2} =
3,53642$, which is different in modulus from the other
eigenvalues.

\section{Remarks}

The results of this article can be easily generalized to the case
of $k$-totally ${\mathcal J}$-sign-symmetric matrices with $k = 3,
\ 4, \ 5, \ \ldots.$

\newpage

\section*{ References}

\medskip

1. {\it F.R. Gantmacher, M.G. Krein.} Oscillation Matrices and
Kernels and Small Vibrations of Mechanical Systems.
--- AMS Bookstore, 2002. --- 310 p.

2. {\it M. Fiedler, V.P. Pt\'{a}k.} Some generalizations of
positive definiteness and monotonicity
// Numerische Math. --- 1966. --- Vol. 9 --- P.
163-172.

3. {\it R.B. Kellog.} On complex eigenvalues of $M$ and $P$
matrices
// Numerische Math. --- 1972. --- Vol. 19 --- P.
170-175.

4. {\it D.M. Kotelyanskii.} On some sufficient criteria for
reality and simpleness of a matrix spectrum
// Mat. Sbornik --- 1955. --- Vol. 36(78), № 1 --- P.
163--168. (Russian)

5. {\it L.M. DeAlba, T.L. Hardy, L. Hogben, A. Wangness.} The
(weakly) sign symmetric $P$-matrix completion problem // Electron
J. Linear algebra. --- 2003. --- Vol. 10 --- P. 257-271.

6. {\it S.M. Fallat, C.R. Johnson.} Multiplicative principal-minor
inequalities for tridiagonal sign-symmetric $P$-matrices.
// Taiwanese Journal of Mathematics. --- 2001. --- Vol. 5, No. 3 --- P.
655--665.

7. {\it D.M. Kotelyanskii.} On a property of sign-symmetric
matrices.
// UMN --- 1953. --- Vol. 8, № 4(56) --- P.
163--167. (Russian)

8. {\it D. Hershkowitz, N. Keller.} Spectral properties of sign
symmetric matrices // Electron J. Linear algebra. --- 2005. ---
Vol. 13 --- P. 90-110.

9. {\it O. Holtz, H. Schneider.} Open problems on GKK
$\tau$--matrices.
//  Linear Algebra and its Applications. --- 2002.
Vol. 345 --- P. 263-267.

10. {\it D. Carlson.} A class of positive stable matrices
// J. Res. Nat. Bur. Stand. --- 1974. --- Vol. 78B --- P. 1--2.

11. {\it Tsoy-Wo Ma.} Classical analysis on normed spaces. ---
World Scientific Publishing, 1995. --- 356 p.

12. {\it I.M. Glazman, Yu.I. Liubich.} Finite-dimensional linear
analysis. - M.: Nauka, 1969.- 476 p. (Russian)

13. {\it J.L. Kelley.} General Topology. --- Birkh\"{a}user, 1975.
--- 298 p.

14. {\it K. Kuratovski.} Topology, I, II, revised 2nd ed. ---
Academic Press, New York, 1966.

15. {\it F. Gantmacher.} The Theory of Matrices. Volume 1, Volume
2.
--- Chelsea. Publ. New York, 1990.

\end{document}